\documentclass[a4paper,12pt,oneside]{amsart}

\def\call#1{\ensuremath{{\mathcal #1}}}
\def\M#1{\ensuremath{\mathbb #1}}
\def\mfk#1{\ensuremath{{\mathfrak{#1}}}}
\def\f{\ensuremath{\varphi}}
\def\e{\ensuremath{\varepsilon}}
\def\wt#1{\ensuremath{\widetilde{#1}}}

\def\wb#1{\ensuremath{\overline{#1}}}

\def\z{{\bf z}}
\def\vol{{\rm vol}}
\def\B{{\rm bar}}

\def\isomk{\ensuremath{{\rm Isom(\M H^k)}}}
\def\isomn{\ensuremath{{\rm Isom(\M H^n)}}}
\def\hKN{\ensuremath{\wb {\M H}^k\times\wb{\M H}^n}}

\def\intbK{\ensuremath{\int_{\partial \M H^k}}}
\def\intbN{\ensuremath{\int_{\partial \M H^n}}}
\def\intK{\ensuremath{\int_{\wb {\M H}^k}}}
\def\intN{\ensuremath{\int_{\wb {\M H}^n}}}
\def\intbKN{\ensuremath{\iint_{\partial \M H^k\times\wb{\M H}^n}}}
\def\intKN{\ensuremath{\iint_{\wb {\M H}^k\times\wb{\M H}^n}}}

\def\Id{{\rm Id}}

\newtheorem{teo}{Theorem}[section]
\newtheorem{prop}[teo]{Proposition}

\newtheorem{cor}[teo]{Corollary}
\newtheorem{lemma}[teo]{Lemma}
\newtheorem{defi}[teo]{Definition}
\newtheorem{remark}[teo]{Remark}
\newtheorem{examp}[teo]{Example}

\newcounter{step}

0

\begin{document}
\title{Constructing equivariant maps for representations}
\author{Stefano Francaviglia}
\thanks{This work was supported by INdAM and the European
  Research Council (MEIF-CT-2005-010975 and MERG-CT-2007-046557.)} 
\address{Dipartimento di Matematica Applicata ``U.Dini'', via
  Buonarroti 1/c, 56127 Pisa, Italy}
\email{s.francaviglia@sns.it}

\begin{abstract}
We show that if
$\Gamma$ is a discrete subgroup of the group of the isometries of $\M
H^k$, and if $\rho$ is a representation of $\Gamma$ into the group of
the isometries of $\M H^n$, then any
$\rho$-equivariant map $F:\M H^k \to \M H^n$ extends to the boundary
in a weak sense in the setting of Borel measures.
As a consequence of this fact, we obtain an extension of a result of
Besson, Courtois and Gallot about the existence of volume
non-increasing, equivariant maps.
Then, we show that the weak
extension we obtain is actually a measurable $\rho$-equivariant
map in the classical sense. We use this fact to obtain measurable
versions of Cannon-Thurston-type results for equivariant Peano
curves. For example, we prove that if $\Gamma$ is of divergence
type and $\rho$ is non-elementary, then there exists a
measurable map $D:\partial \M H^k\to\partial\M H^n$ conjugating the
actions of $\Gamma$ and $\rho(\Gamma)$.
Related applications are discussed.
\end{abstract}

\maketitle
 \tableofcontents

\section{Introduction}
Let $\Gamma$ be a discrete subgroup of \isomk\ and let
$\rho:\Gamma\to\isomn$ be a representation. Then, it is easy to construct a
piecewise smooth map $D:\M H^k\to\M H^n$ which is $\rho$-equivariant,
that is $D(\gamma x)=\rho(\gamma)D(x)$ for all $\gamma\in\Gamma$ and
$x\in\M H^k$, and the problem arises of whether such a map continuously
extends to
the boundaries of the hyperbolic spaces (this is a key step in the
proofs of rigidity results for hyperbolic manifolds, see for
example~\cite{Thu:note,BePe:libro,dun99,fra04,Kla:tesi}). In some cases,
for example if the representation is not discrete, such an extension
is not possible. Moreover, in general it is hard even to construct a
$\rho$-equivariant map between the
boundaries with some regularity properties like continuity and
measurability.
If $\Gamma<{\rm Isom}(\M H^2)$ is a surface group and if
$\rho:\Gamma\to{\rm Isom}(\M H^3)$ is an isomorphism such that $\M
H^3/\rho(\Gamma)$ is an hyperbolic $3$-manifold, under certain assumptions,
Cannon and Thurston~\cite{cathu89} and Minsky~\cite{min94},
proved the existence of
a continuous, $\rho$-equivariant, surjective map $\partial \M
H^2\to\partial \M H^3$, and Soma~\cite{som95} proved that, outside
zero-measure sets, such a map is a homeomorphism
(See also the recent works~\cite{Mjpre06,Mjpre07}.)

\ \\
The starting point of this paper is the existence, for any
$\rho:\Gamma\to\isomn$, of a measurable, $\rho$-equivariant map from
the limit set of $\Gamma$ to the set of probability measures on
$\partial \M H^n$, that can be stated as follows (see
Section~\ref{s_2} for precise definitions.)

\begin{teo}[Existence of developing measures]
  \label{t_3}
Let $\Gamma<\isomk$ be an infinite, non-elementary discrete group.
Let $\rho:\Gamma\to\isomn$
be a representation. Then, a family of developing measures
for $\rho$ exists.
\end{teo}

The fact that we work in the world of measure with the weak-*
topology introduces a lot of compactness, so that we will able to
establish results of existence and convergence for equivariant maps.
The main applications of our technique belong to the framework of
{\em
  Barycentric maps} and the one of {\em Cannon Thurston maps}.

\vskip\baselineskip
\noindent\textsc{Barycentric maps.} First of all, we get a
generalisation of the celebrated Theorem of Besson, Courtois and
Gallot on existence of natural maps.

\begin{teo}[Existence of B-C-G-natural maps]\label{t_1}
  Let $\Gamma<\isomk$ be an infinite discrete group. Let $\rho:\Gamma\to
  \isomn$ be a representation whose image is non-elementary. Then, there
  exists a map $F:\M H^k\to \M H^n$, called {\em natural map},  such that:
  \begin{enumerate}
    \item $F$ is smooth.
    \item $F$ is $\rho$-equivariant, i.e. $F(\gamma
    x)=\rho(\gamma)F(x)$ for all $x\in\M H^k$ and
    $\gamma\in\Gamma$.
  \item\label{t1_p3}  For all $p\geq 3,\  \displaystyle{
{\rm Jac}_pF(x)\leq\left(\frac{\delta(\Gamma)}{p-1}\right)^p}$.
\item\label{t1_p4} If
$||d_xF(u_1)\wedge\dots\wedge d_xF(u_p)||=\displaystyle{
\left(\frac{\delta(\Gamma)}{p-1}\right)^p}$ for an orthonormal $p$-frame
$u_1,\dots,u_p$ at $x\in\M H^k$, then the restriction of $d_xF$ to the
subspace generated by $u_1,\dots,u_p$ is a homothety.
  \end{enumerate}
\end{teo}

Theorem~\ref{t_1} was proved by Besson
Courtois and Gallot in the special case in which $\rho$ is discrete
and faithful and
both $\Gamma$ and $\rho(\Gamma)$ are  convex co-compact
(\cite{bcg99}.)

In our
setting, Theorem~\ref{t_1} will follow directly from Theorem~\ref{t_3}
and a modification of the construction of Besson, Courtois and Gallot (see
\cite{bcg95,bcg96,bcg99}.)

In~\cite{bcg99},
for $\e>0$ the authors construct a smooth
$\rho$-equivariant map
$F_\e:\M H^k\to\M H^n$ such that for all $p\geq 3,\  \displaystyle{
{\rm Jac}_pF_\e(x)\leq\left(\frac{\delta(\Gamma)(1+\e)}{p-1}\right)^p}$.
We call such maps \e-natural maps.
We will see that the natural map we construct is actually the
limit of a sequence of such \e-natural maps.

\begin{teo}\label{p_5}
  In the hypotheses of Theorem~\ref{t_1}, there exists a family
  $\{F_\e\}$ of \e-natural maps (constructed as in~\cite{bcg99})
 and a sequence $\e_i\to 0$ such that
  $F_{\e_i}$ converges to the natural map $F$.
\end{teo}

The proofs of Theorems~\ref{t_3} and~\ref{p_5} both start with a
$\rho$-equivariant map $D:\M H^k\to\M H^n$. While the natural map $F$
does not depend on $D$, the maps $F_\e$'s can be constructed in such a way to
 keep memory of $D$.
More precisely, the
construction of the $\e$-natural maps depends on the choice of a
probability measure on a fundamental domain for the action of
$\Gamma$, and they depends on the restriction of $D$ to the support of
such measure.
This is useful to
study non-compact manifolds.
For example, suppose that $f:M\to N$ is a
proper map between complete non-compact hyperbolic manifolds. If $D$
is the lift of $f$ to the universal covers, then the natural maps
$F_\e$ descend to maps $f_\e:M\to N$, and one can show that, if one
used a suitable measure to construct such maps, then they
are proper.

Such results can be used to prove rigidity results for
representations. The following theorem (whose proof will be sketched
in Section~\ref{s_R} and completely described
in~\cite{frkl06})
 is an example of applications of
Theorems~\ref{t_1} and~\ref{p_5}.

\begin{teo}[Rigidity of representations]\label{t_6.0_rigidity}
Let $M$ be a complete hyperbolic $k$-manifold of finite volume. Let
$\rho:\pi_1(M)\to \isomn$ be an irreducible representation. Then
$\vol(\rho)\leq\vol(M)$ and equality holds if and only if $\rho$ is a
discrete and faithful representation into the group of isometries of a
$k$-dimensional hyperbolic subspace of $\M H^n$.
\end{teo}

\vskip\baselineskip
\noindent\textsc{Measurable Cannon-Thurston maps.}
A Cannon-Thurston map for $\rho$ is a continuous, $\rho$-equivariant
map from the limit set of $\gamma$ to $\partial \M H^n$. It is quite
difficult to show the existence of such maps (and in general there are
obstructions to continuity, and one needs to impose geometric constraints.) In
fact, the existence problem has been solved for $k=2$ and $n=3$ with
geometric hypotheses (see
for instance~\cite{cathu89,Mcm01,min94,Mjpre06,Mjpre07}.)

Theorem~\ref{t_3} is a weak existence result, and its proof
use a kind of weak extension result; namely, the family of developing
measures weakly extends the orbit of a point.
The main result here, is the proof
of a stronger measurable
extension result, which will be the base for the study of
measurable versions of the Cannon-Thurston map (the measures we consider on
$\partial \M H^k$
are the Patterson-Sullivan measures, see  Section~\ref{s_2}.) The
strong extension theorem can be stated in his general form as follows.

\begin{teo}[Existence and uniqueness of measurable extensions]\label{t_4}
Let $\Gamma<\isomk$ be a discrete group and let $\rho:\Gamma\to\isomn$ be a
representation whose image is non-elementary.
Suppose that there exists a family $\{\lambda_z\}_{z\in\partial \M
  H^k}$ of developing measures
for $\rho$  such that for almost all $z$, the measure
$\lambda_z$ is not the sum of two Dirac
deltas with equal weights. Then, the natural map $F$ constructed using
$\{\lambda_z\}$ extends to the conical limit set. More precisely,
there exists a measurable, $\rho$-equivariant map $\wb
F:\wb{\M H}^k\to\wb{\M H}^n$ which agrees with $F$ on $\M H^k$ and
such that, for almost all $\omega$ in the conical limit set,
if $\xi\in\M H^k$ and $\{\gamma_n\xi\}$ is
a sequence conically converging to $\omega$, then
$F(\gamma_n x)\to\wb F(\omega)$ for all $x\in\M H^k$.

Moreover, the map $\wb F$ is unique. More precisely, if $\wb F_1$
and $\wb F_2$ are two measurable,
$\rho$-equivariant maps from $\Lambda(\Gamma)$ to the set of
probability measures on $\partial \M H^n$,
then $\wb F_1$ and $\wb F_2$ are in fact ordinary functions, that is to say,
they map almost every point of $\Lambda(\Gamma)$ to
a Dirac delta concentrated on a point of $\partial\M H^n$. Moreover,
 they agree almost everywhere on $\partial \M H^k$.
\end{teo}

Theorem~\ref{t_4}
 in particular applies to the case of fundamental groups of
hyperbolic
manifolds. Indeed, we will prove that in this case the developing
 measures are almost never the sum of two Dirac
deltas with equal weights. In particular, this gives the following theorem.

\begin{teo}[Existence and uniqueness of Cannon-Thurston maps]\label{c_5.1_1.9}
  Let $\M H^k/\Gamma$ be a complete hyperbolic manifold of finite
  volume and let $\rho:\Gamma\to\isomn$ be a representation whose
  image is not elementary. Then there exists a measurable
  $\rho$-equivariant map $\wb
  F:\partial\M H^k\to\partial\M H^n$
Such a map is the extension of a natural
  map, and two such maps are equal almost everywhere.
\end{teo}

In particular when a classical Cannon-Thurston map exists, it
coincides with the extension provided by Theorem~\ref{c_5.1_1.9}.
An immediate consequence of these facts is the following result.

\begin{teo}[Inverse of Cannon-Thurston maps]\label{t_5.1_7.6}\ \\
Let $\Gamma<\isomk$ and $\Gamma'<\isomn$ be non-elementary discrete
groups such that they diverge respectively  at $\delta(\Gamma)$ and
$\delta(\Gamma')$.
Let $\rho:\Gamma\to\Gamma'$ be an isomorphism.

Then, there exist measurable maps
  $\wb F:\partial \M H^k\to\partial\M H^n$ and  $\wb G:\partial \M
  H^n\to\partial\M H^k$ which are respectively $\rho$ and
  $\rho^{-1}$ equivariant. Moreover, almost everywhere
$$\wb F\circ\wb G={\rm Id}_{\M H^k}\qquad \wb G\circ\wb F={\rm Id}_{\M H^n}.$$
\end{teo}

As noticed above, working with measures helps when one has to deal
with convergence problems. For example, an application of
the above techniques provides an answer to the following question.

\ \\
Suppose that $\rho_i\to\rho$ is a converging sequence of
non-elementary representations. Do the corresponding Cannon-Thurston
maps converge?

\ \\
Miyachi proved~\cite{Miy06} that if $\Gamma$ is a surface-group
without parabolics, and if the injectivity radius is bounded away
from zero along the whole sequence, then the answer is "yes" (in
fact, one has uniform convergence,) and it is conjectured that the
same holds in the case of cusped surfaces with a uniform bound on
the injectivity radius outside the cusps (see for
example~\cite[Conjecture~5.2]{scorzapre}.) In general, as
Example~\ref{ex_7.2_9.1}
 shows, one can construct sequences having no
uniformly converging sub-sequence (the condition on injectivity
radii is violated.) We prove here that in the general case (no
geometric hypotheses, no bounds on dimensions) the answer is "yes"
for the point-wise convergence almost everywhere.

\begin{teo}[Convergence of Cannon-Thurston maps]\label{t_7.0_1.8}
  Let $\Gamma<\isomk$ be a discrete group that diverges at its
  critical exponent $\delta(\Gamma)$. Let $\rho_i:\Gamma\to\isomn$ be
  a sequence of representations with non-elementary images. Suppose
  that $\rho_i$ converges to a representation $\rho$ whose image is
  non-elementary. Let $f_i$ and $f$ be the corresponding measurable
  Cannon-Thurston maps for $\rho_i$ and $\rho$ respectively. Then,
  $f_i$ converges to $f$ almost everywhere with respect to the
  Patterson-Sullivan measures on the limit set of $\Gamma$.
\end{teo}

\ \\
\textsc{Acknowledgements.} I need to thank Jeff Brock, Ben
Klaff, Carlo Mantegazza, Joan Porti, Juan Souto and Xavi Tolsa for
their fundamental help.


\section{Definitions, notation and preliminary facts}\label{s_2}

First of all, we recall some definitions.
\begin{defi}
  A subgroup $\Gamma<\isomk$ is said {\em non-elementary} if any
  $\Gamma$-invariant, non-empty set $A\subset \partial\M H^k$ contains at least
  three points. Otherwise, $\Gamma$ is said {\rm elementary}.
\end{defi}

\begin{defi}
Let $F:\M H^k\to\M H^n$ be a smooth map. The $p$-Jacobian ${\rm
  Jac}_pF$ of $F$ is defined by $${\rm Jac}_pF(x)=\sup
  ||d_xF(u_1)\wedge,\dots,\wedge d_xF(u_p)||_{\M H^n},$$
where $\{u_i\}_{i=1}^p$ varies on the set of orthonormal $p$-frames at
  $x\in\M H^k$.
\end{defi}

\begin{defi}
  Let $(X,g)$ be a complete Riemannian manifold and let $\Gamma$ be a
  group of isometries of $X$. We denote by $\delta(\Gamma)$ the
  critical exponent of the Poincar\'e series of $\Gamma$, that is:
$$\delta(\Gamma)=\inf\{s>0 : \sum_{\gamma\in\Gamma} e^{-sd(x,\gamma
  x)}<+\infty\}$$
where $d(\cdot,\cdot)$ denotes the distance induced by $g$ on $X$ and
  $x$ is a point of $X$.
\end{defi}

It is readily checked that $\delta(\Gamma)$ does not depend on $x$.
We notice that, when $X=\M H^k$,
 the critical exponent is the Hausdorff dimension of the conical limit set of
 $\Gamma$ (\cite{BiJo97}). Moreover, $\delta(\Gamma)$ can be computed by
$$\delta(\Gamma)=\lim_{R\to\infty}\frac{1}{R}\log(\#\{\gamma\in\Gamma:\
d(\gamma O,O)<R\}).$$
We refer to~\cite{nic:libro,sul84,yue96annals} for details.
\begin{defi}
 Let $(X,g)$ be a complete Riemannian manifold and let $\Gamma$ be a
  group of isometries of $X$. We say that $\Gamma$ {\em diverges} at
$\delta(\Gamma)$ if
 $$\lim_{s\to\delta(\Gamma)^+}
\sum_{\gamma\in\Gamma} e^{-sd(x,\gamma x)}=+\infty.$$
\end{defi}

The following lemma is an immediate consequence of
\cite[Theorem 1.6.1]{nic:libro}, \cite[Theorem 1.6.3]{nic:libro} and \cite[Corollary 3.4.5]{nic:libro}.
\begin{lemma}\label{l_nic}
Let $\Gamma$ be a infinite, non-elementary discrete  subgroup of $\isomk$, then
$$0<\delta(\Gamma)\leq k-1.$$
Moreover, if $\M H^k/\Gamma$ has finite volume, then
$\delta(\Gamma)=k-1$
and $\Gamma$ diverges at $\delta(\Gamma)$.\end{lemma}

\vskip\baselineskip
\noindent
\textsc{Notation.} For the rest of the paper we fix the following
notation. $\Gamma$ will be an infinite discrete subgroup of \isomk, and
$$\rho:\Gamma\to\isomn$$ will be a representation.
We fix base-points in $\M H^k$ and $\M H^n$, both denoted by $O$.
We denote by $B_K$ and $B_N$ the Busemann
functions, respectively of $\M H^k$ and $\M H^n$, normalised at
$O$. Namely, for $x\in\M H^k$ and $\theta\in\partial\M H^k$ (resp. $\M H^n$
and $\partial \M H^n$) we set
$$B_K(x,\theta)=\lim_{t\to\infty} \big(d_{\M
  H^k}\big(x,\gamma_\theta(t)\big)-t\big),$$
where $\gamma_\theta$ is the geodesic ray from $O$ to $\theta$,
  parametrised by arc length.
We denote by $\pi_K$ (respectively $\pi_N$) the projection of $\wb {\M
  H}^k\times \wb{\M H}^n$ to $\wb{\M H}^k$ (resp. $\wb{\M H}^n$):
$$\pi_K:\wb {\M H}^k\times \wb{\M H}^n\to \wb{\M H}^k.$$
Finally, we fix a continuous piecewise smooth  $\rho$-equivariant map
$$D:\M H^k\to \M H^n.$$
We notice that such a map can be easily constructed by
triangulating a fundamental domain for $\Gamma$ and then arguing
by induction on the $i$-skeleta.

\vskip\baselineskip
\noindent
\textsc{Patterson-Sullivan measures.}
A fundamental tool for our purpose is the family of
Patterson-Sullivan measures. We recall the main results we need, which
we summarise in Theorem~\ref{t_psm}, referring
to~\cite{bcg99}, \cite[Chap. 3, 4]{nic:libro} and~\cite{yue96annals,yue96} for proofs and details.

\begin{teo}\label{t_psm} Let $\Gamma$ be an infinite, non-elementary discrete
  subgroup of \isomk\ with critical exponent $\delta(\Gamma)$.
  For all $x\in\M H^k$ there exists a positive Borel measure $\mu_x$
  of finite, non-zero mass such
  that, for all $x,y\in\M H^k$ and $\gamma\in\Gamma$:
  \begin{enumerate}
  \item\label{psm_1}
The measures $\mu_x$ and $\mu_y$ are in the same density class
  of measures and are concentrated on $\partial \M H^k$.
\item\label{psm_2} The measure $\mu_x$ satisfies
  $$d\mu_x(\theta)=e^{-\delta(\Gamma)B_K(x,\theta)} d\mu_O(\theta)$$
where $\theta\in\partial \M H^k$.
\item\label{psm_3} The measures $\mu_x$ are $\Gamma$-equivariant, that is
$$\mu_{\gamma x}=\gamma_*\mu_x.$$
  \end{enumerate}
\end{teo}

\proof
We only sketch the proof, recalling  the construction of
the Patterson-Sullivan measures because we will explicitly
use it in the following.

For all $s>\delta(\Gamma)$ let
$$c(s)=\sum_{\gamma\in\Gamma}e^{-sd(O,\gamma O)}$$
where $d(\cdot,\cdot)$ denotes the hyperbolic distance of $\M H^k$.
 For simplicity, here we stick to the case  that $\Gamma$ diverges at
 $\delta(\Gamma)$
 (this happens for example if $\Gamma$ is geometrically finite,
 see~\cite[p. 87]{nic:libro}).
For any $x\in\M H^k$ and $s>\delta(\Gamma)$, define
$$\mu_x^s=\frac{1}{c(s)}\sum_{\gamma\in\Gamma}e^{-sd(x,\gamma
  O)}\delta_{\gamma O}$$
where $\delta_{\gamma O}$ is the Dirac measure concentrated on $\gamma
  O$. For $s>\delta(\Gamma)$, $\mu_x^s$ is a well-defined positive
  Borel measure on $\M H^k\subset \wb{\M H}^k$. It can be shown
  that for $s\to\delta(\Gamma)^+$, the
  measures $\mu_x^s$ weakly converge to a positive Borel measure
  $\mu_x$ on $\wb{\M H}^k$. More precisely, for all $x\in\M H^k$ and
  $\f\in C(\wb{\M H}^k)$
$$\intK \f\, d\mu_x^s\to\intK \f\, d\mu_x.$$
Moreover, the fact that $\lim_{s\to\delta(\Gamma)^+}c(s)=+\infty$
easily implies that $\mu_x$ is concentrated on the boundary (in fact
$\mu_x$ is concentrated on the limit set of $\Gamma$).
\qed

\vskip\baselineskip
\noindent
\textsc{Barycentre of a measure.}
We recall now the definition and the main properties of the barycentre
of a measure. We refer the reader to~\cite{DoEa86} and~\cite{bcg95}
for complete proofs and details.

Let $\beta$ be a positive Borel measure on $\partial \M H^n$ of finite
mass. Define
the function $\call B_\beta:\M H^n\to\M R$ by
$$\call B_\beta: y\mapsto \intbN B_N(y,\theta)\,d\beta(\theta).$$

Since we are working in the hyperbolic space, the Busemann functions
are convex. Thus, its $\beta$-average $\call B_\beta$ is strictly convex,
provided that $\beta$ is not the sum of two deltas. Moreover, one
can show that $\call B_\beta(y)\to\infty$ as $y$ approaches $\partial \M
H^n$. It follows that $\call B$ has a unique minimum in $\M H^n$.
\begin{defi}
  For any positive Borel measure $\beta$ on $\partial \M H^n$ of
  finite mass which is not concentrated on two points, we define the
  {\em barycentre} $\B(\beta)$ of $\beta$ as the unique minimum point
  of the function $y\mapsto \intbN B_N(y,\theta)d\beta(\theta)$.
\end{defi}

We refer to~\cite{DoEa86},~\cite{bcg95} and~\cite{bcg96} for a proof
of the following lemma.

\begin{lemma}\label{l_bar}
 The barycentre of a measure $\beta$ satisfies the following
 properties:
 \begin{enumerate}
 \item\label{eq_bar1}
The barycentre is characterised by the equation $$\displaystyle{\intbN
 {dB_N}_{(\B(\beta),\theta)}(\cdot)\,d\beta(\theta)=0.}$$
\item\label{eq_bar2} The barycentre is \isomn-equivariant, that is, for any
  $g\in\isomn$ $$\B(g_*\beta)=g(\B(\beta)).$$
\item\label{eq_bar3} The barycentre is continuous respect the weak convergence of
  measures. That is, if $\beta_i\rightharpoonup\beta$, then
  $\B(\beta_i)\to\B(\beta)$.
 \end{enumerate}
\end{lemma}
The first property follows from the definition
after  differentiating the function $\call B_\beta$.
The equivariance follows from
the properties of the Busemann functions. The continuity can be easily
proved using that, if $\beta_i\rightharpoonup\beta$, then $\call
B_{\beta_i}$ and $d\call B_{\beta_i}$ point-wise converge to $\call
B_\beta$ and $d\call B_\beta$ respectively.

\begin{remark}
If $\beta=a\delta_{\theta_1}+b\delta_{\theta_2}$, with $0<a<b$, then it
can be checked that the minimum of $\call B_\beta$ is the point
$\theta_2\in\partial \M H^n$.
Thus, one can define the barycentre of a measure $\beta$
whenever $\beta$ is not the sum of two deltas with the same weights.
Note that, since the barycentre of a measure concentrated on two points
belongs to $\partial \M H^n$, equation~\ref{eq_bar1} of
Lemma~\ref{l_bar} makes no sense for such measures.

\end{remark}
\vskip\baselineskip
\noindent
\textsc{Developing measures.}
 We now introduce the notion of {\em family of developing measures}
 for $\rho$, which extends the one of $\rho$-equivariant map.
We recall that
 $\{\mu_x\}$ is the family of Patterson-Sullivan measures.

 \begin{defi}\label{d_devm}
   A family of {\em developing measures} for $\rho$ is a set
   $\{\lambda_z\}_{z\in\partial \M H^k}$ of positive Borel measures on
   $\partial \M H^n$, of finite mass, and such that:
   \begin{enumerate}
   \item\label{devm_1}
 The measures $\lambda_z$'s are $\rho$-equivariant, that is, for
   $\mu_O$-almost all $z$ and all $\gamma\in\Gamma$
$$\lambda_{\gamma z}=\rho(\gamma)_*\lambda_z$$
\item\label{devm_2} For any $\f\in C(\partial \M H^n)$, the function
$$z\mapsto \intbN \f(\theta)\,d\lambda_z(\theta)$$ is $\mu_O$-integrable
(whence, by points~\ref{psm_1} and \ref{psm_2} of Theorem~\ref{t_psm},
it is $\mu_x$-integrable for all $x$).
\item\label{devm_3} The function $z\mapsto||\lambda_z||$ belongs to
  $L^\infty(\partial \M H^k,\mu_O)$.
\item\label{devm_4} For $\mu_O$-almost $z\in\partial \M H^k$,
$||\lambda_z||>0.$
   \end{enumerate}
 \end{defi}

As an example, consider a $\mu_O$-measurable  $\rho$-equivariant map
$\wb D:\partial \M H^k\to\partial \M H^n$. Then the family
$\{\lambda_z=\delta_{\wb D(z)}\}$, where $\delta_{\wb D(z)}$ is the
Dirac measure, is a family of developing measures for $\rho$. In this
sense the notion of developing measures extends the one of equivariant
map.

\vskip\baselineskip
\noindent
\textsc{Convolutions of measures.}
Let $X,Y$ be topological spaces and let $\mu$ be a Borel measure on
$X$. Let $\{\alpha_x\}_{x\in X}$ be a family
of Borel measures on $Y$ such that for each
$\f\in C_0(Y)$ the function $x\mapsto\int_Y\f
\,d\alpha_x$ is $\mu$-integrable.
The convolution $\mu*\{\alpha_x\}$ is the Borel measure on $Y$ defined by
$$\int_Y \f(y)\,d(\mu*\{\alpha_x\})=\int_X\left(\int_Y\f(y)\,d\alpha_x(y)\right)\,d\mu(x)$$
for any $\f\in C_0(Y)$. Similarly, we define the product
$\mu\times\{\alpha_x\}$ on $X\times Y$ by
$$\int_Y \f(x,y)\,d(\mu\times\{\alpha_x\})=\int_X\left(\int_Y\f(x,y)\,d\alpha_x(y)\right)\,d\mu(x).$$
The measure $\mu*\{\alpha_x\}$ is the $\mu$-average of the
$\alpha_x$'s. Moreover, if $\pi:X\times Y\to Y$ is the projection,
then $\mu*\{\alpha_x\}=\pi_*(\mu\times\{\alpha_x\})$.

We say that a sequence of measures $\{\mu_i\}$ weakly converges to
$\mu$ if, for any continuous function $f$
with compact support,
$\int\!\!f\,d\mu_i\to\int\!\!f\,d\mu$.
The proof of following lemmas are left to the reader.
\begin{lemma}\label{l2.0_2.9}
Suppose that $\{\mu_i\}$ is a sequence of measures on $X$,
weakly converging to $\mu$.
If for each $\f\in C_0(Y)$ the function $x\mapsto\int_Y\f
\,d\alpha_x$ belongs to $C_0(X)$, then the sequence
$\mu_i*\{\alpha_x\}$ weakly converges to $\mu*\{\alpha_x\}$.
\end{lemma}

\begin{lemma}\label{l2.0_2.10}
Let $Z$ be a topological space and let $\{\nu_y\}_{y\in Y}$ be a
family of Borel measures on $Z$ such that for all $\psi\in C_0(z)$
the function $y\mapsto\int_Z\psi
\,d\nu_y$ is $\alpha_x$-integrable for $\mu$-almost all $x$ and
$x\mapsto \int_Y\int_Z\psi\,d\nu_yd\alpha_x$ is $\mu$-measurable. Then
$$\mu*\{\alpha_x*\{\nu_y\}\}=(\mu*\{\alpha_x\})*\{\nu_y\}.$$
\end{lemma}
By Lemma~\ref{l2.0_2.10} we can omit the parentheses and write
$\mu*\{\alpha_x\}*\{\nu_y\}$.

\vskip\baselineskip
\noindent
\textsc{Some remarks on the hypotheses of Theorems~\ref{t_3},
  \ref{t_1} and \ref{t_2}.}
  First of all, we notice that the hypotheses of such theorems can
  be relaxed by replacing
  the spaces  $\M H^k$ and $\M H^n$ with Riemannian
  manifolds with suitable bounds on the curvatures. We refer the
  reader to~\cite{bcg99} for further details on that direction.

Moreover, even if the hypothesis that the image of $\rho$ is
non-elementary is crucially used in Corollary~\ref{c5_3.4}, it is not
strictly necessary.
More precisely, suppose that $\rho(\Gamma)$ is elementary.
Then there exists a $\rho(\Gamma)$-invariant set
  $A\subset\partial \M H^n$ with either one or two points.

In the
  latter case, there exists a whole geodesic $c$ in $\M H^n$ which is
  $\rho(\Gamma)$-invariant and it is easy to see that a natural map
  whose image is contained in $c$ exists. Therefore, Theorem~\ref{t_1}
  is true in this case.

In the former case, we can suppose that, in the half-space model $\M
H^{n-1}\times\M R^+$ of $\M H^n$, the point $\infty$
is fixed by $\rho(\Gamma)$. Then, one can easily construct a
$\rho$-equivariant map $D$ whose image is contained in the horosphere
$\{(z,1): z\in\M R^{n-1}\}$. Thus, if ${\rm Jac}_pD(x)$ is bounded,
then Theorem~\ref{t_1} is proved by raising $D$ to a sufficiently
high horosphere. It follows that, if $\M H^k/\Gamma$ is compact, or
simply if $\M H^k/\Gamma$
retracts to a compact set (and this is the case if for example $\Gamma$ is
geometrically finite), then Theorem~\ref{t_1} is true even if
$\rho(\Gamma)$ has a fixed point.


\section{Construction of B-C-G natural maps}\label{s_3}

In this section we prove the following Theorem, which, together with
 Theorem~\ref{t_3} (proved in Section~\ref{s_4},)  directly implies
 Theorem~\ref{t_1}. For this section we keep the notation fixed in
 Section~\ref{s_2}.

\begin{teo}
  \label{t_2}
Let $\Gamma<\isomk$ be an infinite discrete group. Let $\rho:\Gamma\to\isomn$
be a representation whose image is non-elementary. If there exists a
family $\{\lambda_z\}$ of developing measures
for $\rho$, then a natural map exists.
\end{teo}
\proof
Fix a family
$$\{\lambda_z\}_{z\in\partial \M H^k}$$ of developing measures for
$\rho$.
The idea of the proof is to use the developing measures  to
push-forward the Patterson-Sullivan measures $\mu_x$'s to measures
$\beta_x$'s on
$\partial \M H^n$,
and define the natural map by
$$
x\mapsto \mu_x\mapsto \beta_x \mapsto \B(\beta_x).
$$
Then, the properties of the natural map will follow as
in~\cite{bcg99}.

The push-forward the measures $\mu_x$'s is defined as follows.
For each $x\in\M H^k$, define $\beta_x$
as the positive Borel measure on $\partial \M H^n$ given by
$\beta_x=\mu_x*\{\lambda_z\}$. Namely,
for all $\f\in C(\partial \M H^n)$
\begin{eqnarray*}
\intbN
\f(y)\,d\beta_x(y)&=&
\intbK\left(\intbN\f(y)\,d\lambda_z(y)\right)\,d\mu_x(z)\\
&=&\intbK\left(\intbN\f(y)\,d\lambda_z(y)\right)
e^{-\delta(\Gamma)B_K(x,z)}\,d\mu_O(z).
\end{eqnarray*}

Note that the measure $\beta_x$ is well-defined and has finite mass because of
conditions~$(\ref{devm_2})$ and $(\ref{devm_3})$ of
Definition~\ref{d_devm}.
Moreover, since $\beta_x=\pi_{N*}(\mu_x\times\{\lambda_z\})$, if
the family $\{\lambda_z\}$ is of the form $\{\delta_{\wb D(z)}\}$ for a
$\mu_O$-measurable function $\wb D:\partial \M H^k\to\partial\M H^n$,
then $\beta_x=\wb D_*\mu_x$.

\begin{lemma}\label{l_bequi}
  The family of measures $\{\beta_x\}_{x\in\M H^k}$ is
  $\rho$-equivariant, that, is for all $x\in\M H^k$ and
  $\gamma\in\Gamma$
$$\beta_{\gamma x}=\rho(\gamma)_*\beta_x.$$
\end{lemma}
\proof
For any $\f\in C(\partial \M H^n)$
\begin{eqnarray*}
 && \intbN\f(y)\,d\beta_{\gamma
  x}(y)=\intbK\left(\intbN\f(y)\,d\lambda_z(y)\right)\,d\mu_{\gamma
  x}(z)\\
&=&\intbK\left(\intbN\f(y)\,d\lambda_{\gamma z}(y)\right)\, d\mu_x(z)\\
&=&\intbK\left(\intbN \f(\rho(\gamma)y)\,
  d\lambda_z(y)\right)\,d\mu_x(z)=\intbN\f(y)\,d(\rho(\gamma)_*\beta_x)(y).
\end{eqnarray*}
\qed
\begin{lemma}\label{l_bden}
  For all $x,y\in \M H^k$, the measures $\beta_x$
  and $\beta_y$ are in the same density class of measures.
\end{lemma}
\proof We have to show that for all positive functions $\f\in C(\partial \M
H^n)$ we have
$$\intbN\f\,d\beta_x=0\ \Longleftrightarrow\ \intbN\f\,d\beta_y=0.$$

This follows from the fact that $\mu_x$ and $\mu_y$ are in the same
density class.
Indeed, if $\Phi$ denotes the function $z\mapsto\intbN\f\,d\lambda_z$,
since the developing measures are positive, $\Phi$ is positive,
and
$$\intbN\f\,d\beta_x=\intbN\Phi\,d\mu_x=0
\ \Longleftrightarrow\
0=\intbN\Phi\,d\mu_y=\intbN\f\,d\beta_y.$$
\qed

\begin{lemma}\label{l_b>}
  For all $x\in\M H^k$, $||\beta_x||>0$.
\end{lemma}
\proof
It follows from condition~$(\ref{devm_4})$ of Definition~\ref{d_devm}
and from the fact that $||\mu_x||>0$.\qed
\begin{cor}\label{c5_3.4}
  For all $x\in\M H^k$, the measure $\beta_x$ is not concentrated on
  two points.
\end{cor}
\proof By Lemma~\ref{l_b>}, $\beta_x$ is not the zero-measure.
Suppose that $\beta_x$ has an atom of positive weight at
$y_0\in\partial \M H^n$. By Lemmas~\ref{l_bequi} and~\ref{l_bden},
$\beta_x$ has an
atom of positive weight at each point of the $\rho(\Gamma)$-orbit of
$y_0$, which contains most than two points because $\rho(\Gamma)$ is
non-elementary. \qed

It follows that for all $x\in\M H^k$, the barycentre of the measure
$\beta_x$ is well-defined and belongs to $\M H^n$. We define the natural
map $F:\M H^k\to \M H^n$ by
$$F(x)=\B(\beta_x).$$

By condition~$(\ref{eq_bar1})$ of Lemma~\ref{l_bar}, the
natural map is characterised by the implicit equation
\begin{eqnarray}
\label{eq_impl}
G(x,F(x))=0
\end{eqnarray}
where
\begin{eqnarray*}
G(x,\xi)&=&\intbN{dB_N}_{(\xi,\theta)}(\cdot)\,
  d\beta_x(\theta)\\
&=&\intbK\intbN{dB_N}_{(\xi,\theta)}(\cdot)\,
  d\lambda_z(\theta)e^{-\delta(\Gamma)B_K(x,z)}\, d\mu_O(z).
\end{eqnarray*}
The function $G$ is smooth because the Busemann functions $B_K$ and
$B_N$ are smooth. Then, by the implicit function theorem, we get that $F$
is smooth. Moreover, by Lemma~\ref{l_bequi} and
claim~$(\ref{eq_bar2})$ of Lemma~\ref{l_bar}, it follows that $F$ is
$\rho$-equivariant.

By differentiating equation~$(\ref{eq_impl})$ we get that  for all $u\in
T_x\M H^k$ and $v\in T_{F(x)}\M H^n$
\begin{eqnarray}
  \label{eq_bcg_2.3}
  \begin{array}{c}
\displaystyle{\intbK\intbN Dd{B_N}|_{(F(x),\theta)}(dF_x(u),v)\,
  d\lambda_z\,d\mu_x=}\\\ \\
\displaystyle{
=\delta(\Gamma)\intbK\intbN d{B_N}_{(F(x),\theta)}(v)d{B_K}_{(x,z)}(u)
  \,d\lambda_z\,d\mu_x}.
  \end{array}
\end{eqnarray}

Equation~$(\ref{eq_bcg_2.3})$ is the analogous of equation~$(2.3)$
of~\cite{bcg99}. The proof of the properties of the natural map now
goes exactly as in~\cite[p. 152-154]{bcg99},
and the proof of Theorem~\ref{t_2} is complete.
\qed

An immediate corollary of Theorem~\ref{t_1} is the following fact.
\begin{cor}\label{t_6}
  Let $X,Y$ be two compact hyperbolic manifold of (possibly
  different) dimension at least tree. Then in each homotopy class of
  maps $f:Y\to X$ there exists a smooth map $F:Y\to X$ such that
  $|Jac(F)|\leq1$. Moreover, if $|JacF(y)|=1$, then $d_yF$ is an isometry.
\end{cor}
\proof By Theorem~\ref{t_1}, any
$f:Y\to X$ is homotopic to the
natural map corresponding to the representation $f_*$. Such a map has
the requested properties because, by Lemma~\ref{l_nic}, we have
$\delta(\pi_1(Y))\leq{\rm dim}(Y)$.\qed

As Theorem~\ref{t_1}, Corollary~\ref{t_6} was proved by Besson,
Courtois and Gallot in the special case in which $\rho$ is discrete
and faithful and
both $\Gamma$ and $\rho(\Gamma)$ are  convex co-compact
(\cite{bcg99}.)


\section{Weak extension of equivariant maps: existence of developing measures}\label{s_4}

In this section we prove Theorem~\ref{t_3}. In particular, to obtain a
family of developing measures, we show that any equivariant map weakly
extend to the boundary in the setting of Borel measures.
We keep here the notation
fixed in Section~\ref{s_2}, included the one of Theorem~\ref{t_psm} for
the Patterson-Sullivan measures.

The rough idea is the following.
Consider the graphic $G(D)$ of $D$ as a subspace of
  $\wb{\M H}^k\times\wb{\M H}^n$. Then, for each measure $\mu$ on $\M
  H^k$ we can consider the graphic measure $\eta$ on $\wb{\M
  H}^k\times\wb{\M H}^n$, that is, the only measure which is
  concentrated on $G(D)$ and whose push-forward $\pi_{K*}\eta$ is
  $\mu$. The sequence $\{\eta_x^s\}$ of the
  graphic measures corresponding to the measures $\{\mu_x^s\}$,
weakly converges to a measure $\eta_x$ concentrated
  on $\partial \M H^k\times \wb{\M H}^n$. Then, by disintegrating the
  measure $\eta_x$ we obtain a family
  $\{\alpha_z\}_{z\in\partial \M H^k}$ of measures on $\wb{\M H}^n$.
  By making the convolution with the family of visual measures of
  $\partial \M H^n$, we
  get a family of developing measures.

\ \\
For each $x\in\M H^k$ and $s>\delta(\Gamma)$ we define a positive
Borel measure of finite mass on \hKN\ as follows:
\begin{eqnarray*}
  \eta_x^s=\frac{1}{c(s)}\sum_{\gamma\in\Gamma}e^{-sd(x,\gamma
  O)}\delta_{(\gamma O,D(\gamma
  O))}=\frac{1}{c(s)}\sum_{\gamma\in\Gamma}e^{-sd(x,\gamma
  O)}\delta_{(\gamma O,\rho(\gamma)D(O))}
\end{eqnarray*}
where $\delta_{(x,y)}$ denotes the Dirac measure concentrated on
$(x,y)\in\hKN$.  The measures $\{\eta_x^s\}$ are the  graphic measures
associated to $\{\mu_x^s\}$ and $D$, and they are concentrated on $\M
H^k\times\M H^n$. Note that $\eta_x^s=\mu_x^s\times\{\delta_{D(z)}\}$
(see Section~\ref{s_2}.)

The group $\isomk\times\isomn$ acts on \hKN\ by
$$(g_1,g_2)(x,y)=(g_1x,g_2y).$$

\begin{lemma}
  For all $s>\delta(\Gamma)$ the family $\{\eta_x^s\}$ is
  $\rho$-equivariant, that is, for all $x\in \M H^k$ and
  $\psi\in\Gamma$
$$\eta_{\psi x}^s=(\psi,\rho(\psi))_*\eta_x^s.$$
\end{lemma}
\proof
\begin{eqnarray*}
  \eta_{\psi x}^s&=&\frac{1}{c(s)}\sum_{\gamma\in\Gamma}
e^{-sd(\psi x,\gamma O)}\delta_{(\gamma O,D(\gamma O))}
=\frac{1}{c(s)}\sum_{\gamma\in\Gamma}
e^{-sd(x,\psi^{-1}\gamma O)}\delta_{(\gamma O,D(\gamma O))}\\
&=&\frac{1}{c(s)}\sum_{\gamma\in\Gamma}
e^{-sd(x,\gamma O)}\delta_{(\psi \gamma O,D(\psi\gamma O))}=
\frac{1}{c(s)}\sum_{\gamma\in\Gamma}
e^{-sd(x,\gamma O)}\delta_{(\psi\gamma O,\rho(\psi)(D(\gamma O)))}\\
&=&(\psi,\rho(\psi))_*\eta_x^s.
\end{eqnarray*}
\qed

Now, focus on the point $O$ and consider the family $\{\eta_O^s\}_{s>
  \delta(\Gamma)}$. Since 
$||\eta_O^s||=1$ for all $s$, there exists
a sequence $s_i\to\delta(\Gamma)^+$ such that $\eta_O^{s_i}$ weakly
converges to a positive Borel measure $\eta_O$ of finite mass on
$\wb{\M H}^k\times\wb{\M H}^n$
$$\eta_O^{s_i}\rightharpoonup\eta_O.$$
The sequence $\{s_i\}$ depends on the point $O$, but we will show later
that the same sequence works for any point $x$ 
(see Theorem~\ref{t_adisint} below.)
 Moreover, by definition, for all $s>\delta(\Gamma)$ and
  $x\in\M H^k$, the measure $\mu_x^s$ is the $\pi_K$-push-forward of
  $\eta_x^s$. Then, weak continuity of the push-forward implies
$$\mu_O=(\pi_{K})_*\eta_O.$$

In particular, since the support of the Patterson-Sullivan measures is
the limit set 
$\Lambda(\Gamma)$ of $\Gamma$, the measure $\eta_O$ is concentrated on
$\Lambda(\Gamma)\times\wb{\M H}^n$. A priori, $\eta_O$ is non
supported on a graphic; nonetheless, as it is a positive Borel measure,
we can use the theorem of disintegration of measures (\cite[Theorem
2.28]{afp:libro}, see also
~\cite{DeMe:libro},) which asserts that
 for any positive Borel measure $\eta$ on \hKN,
there exists a family of positive Borel measures
$\{\alpha_z^\eta\}_{z\in \wb{\M
    H}^k}$ such that, if $\mu=(\pi_K)_*\eta$, then
$\eta=\mu\times\{\alpha_z^\eta\}$. Thus,
for all $\f\in
C(\hKN)$
$$\intKN
\f\,d\eta=\intK\left(\intN\f(x,y)\,d\alpha_z^\eta(y)\right)\,
d\mu(z),$$
and for $\mu$-almost $z$, the measure $\alpha_z^\eta$ is a
probability measure.
We say that the family $\{\alpha_z^\eta\}$
{\em disintegrates} the measure $\eta$
(compare disintegration with
property~$(\ref{devm_2})$ of Definition~\ref{d_devm}).

Moreover, the measures $\alpha_z^\eta$'s are characterised
by the following property. For all $z\in\wb{\M H}^k$, let
$\{U_j(z)\subset\wb{\M H}^k\}_{j\in\M N}$
be a sequence of nested balls centred
  at $z$  such that $\cap_j U_j(z)=z$. For $j\in\M N$,
let $\psi_j^{z,\eta}:\wb{\M
    H}^k\to \M R$ be the following function (defined for $\mu$-almost
  all $z$) 
$$\psi_j^{z,\eta}=\frac{\chi_{U_j(z)}}{\mu(U_j(z))}$$
where $\chi_A$ denotes the characteristic function of the set
$A$ (note that $\psi_j^{z,\eta}\to\delta_z$ for $\mu$-almost all $z$).
Then, for $\mu$-almost all $z$ and for all $\f\in C(\wb{\M H}^n)$
$$\intN\f(y)\,d\alpha_z^\eta(y)
=\lim_{j\to\infty}\intKN\psi_j^{z,\eta}(x)\f(y)\,d\eta(x,y).$$

From now on, the set $\{\alpha_z\}_{z\in\wb{\M H}^k}$ will denote the
family of measures that disintegrates $\eta_O$, where we set
$\alpha_z=0$ for $z\in\M H^k$. It seems worth mentioning that the
choice of $\alpha_z$ for $z$ 
in the interior of $\M H^k$ is not relevant for our purposes.
Indeed, since any $\mu_x$ is concentrated
on the limit set of $\Gamma$, which is contained in the boundary at
infinity of $\M H^k$, we have $\mu_x(\M{H}^k)=0$. In particular, the
value of $\alpha_z$  for $z$ in the interior of $\M H^k$ does not
affect the product measure $\mu_x\times\{\alpha_z\}$.

\begin{defi}
  For each $x\in\M H^k$ we set $\eta_x=\mu_x\times\{\alpha_z\}$.
\end{defi}

Now, we have two problems. First, while it is true that any
$\{\eta_x^s\}$ has a converging sub-sequence, it is not clear a priori
that the same sub-sequence works for all $x$. Second, by definition, the family
$\{\alpha_z\}$ disintegrates the limit measure $\eta_O$, and there is no
reasons for such a family to disintegrate the limits of $\eta_x^s$ for
all $x$. In other words, it is not clear a priori that $\eta_x^s$
converges to $\eta_x$. The following theorem settles both questions.
This is the core of the proof of Theorem~\ref{t_3}, as it implies that
the family $\{\alpha_z\}$ is equivariant (see Lemma~\ref{l4_4.3} below.)

\begin{teo}\label{t_adisint}
For all $x\in\M H^k$, the sequence $\{\eta_x^{s_i}\}$ weakly converges to
$\eta_x$. Where $\{s_i\}$ is the sequence of indices such that 
$\eta_O^{s_i}\rightharpoonup\eta_O.$
\end{teo}
\proof We fix $x\in\M H^k$.
Since the measures $\eta_x^{s_i}$ are bounded in norm, up to
pass to sub-sequences, we can suppose that $\{\eta_x^{s_i}\}$ weakly
converges to a measure $\wt\eta_x$. We prove now that any such possible
limit $\wt\eta_x$ coincides with $\eta_x$, and this will prove the thesis.
Let $\{\wt\alpha_z\}$ be the family that
disintegrates $\wt\eta_x$. By weak continuity of push forward
$$\pi_{K*}(\wt\eta_x)=\mu_x=\pi_{K*}(\eta_x).$$
Therefore, it is sufficient to show that $\wt\alpha_z=\alpha_z$
for $\mu_O$-almost all $z$ (recall that since $\mu_x$ and $\mu_O$ are
in the same density-class, the notions of $\mu_O$-negligible set and
$\mu_x$-negligible set coincide.)
For any positive function $\f\in C(\wb{\M H}^n)$ and
$z\in\partial \M H^k$
\begin{eqnarray*}
& &  \intN\f(y)\,d\wt\alpha_z(y)
=\lim_{j\to\infty}\intKN\psi_j^{z,\wt\eta_x}(\xi)\f(y)\,d\wt\eta_x(\xi,y)\\
&=&\lim_{j\to\infty}\lim_{s_i\to\delta(\Gamma)^+}\intKN\psi_j^{z,\wt\eta_x}(\xi)\f(y)\,d\eta_x^{s_i}(\xi,y)\\
&=&\lim_{j\to\infty}\lim_{s_i\to\delta(\Gamma)^+}\intK\psi_j^{z,\wt\eta_x}(\xi)\f(D(\xi))\,d\mu_x^{s_i}(\xi)\\
&=&\lim_{j\to\infty}\lim_{s_i\to\delta(\Gamma)^+}
\int_{U_j(z)}\f(D(\xi))\frac{1}{\mu_x(U_j(z))}\,d\mu_x^{s_i}(\xi)\\
&=&\lim_{j\to\infty}\lim_{s_i\to\delta(\Gamma)^+}
\int_{U_j(z)}\f(D(\xi))\frac{\mu_O(U_j(z))e^{-\delta(\Gamma)B_K(x,z)}}{\mu_x(U_j(z))\mu_O(U_j(z))e^{-\delta(\Gamma)B_K(x,z)}}\,d\mu_x^{s_i}(\xi)
\end{eqnarray*}
whence, using the definition of $\mu_x^{s_i}$, and setting
$$A_j(z)=\frac{\mu_O(U_j(z))
e^{-\delta(\Gamma)B_K(x,z)}}
{\mu_x(U_j(z))},$$
we get
\begin{eqnarray}\label{eq_6}
& & \intN\f(y)\,d\wt\alpha_z(y)\\
&=&
\lim_{j\to\infty}A_j(z)\lim_{s_i\to\delta(\Gamma)^+}
\int_{U_j(z)}\frac{\f(D(\xi))}{\mu_O(U_j(z))
e^{-\delta(\Gamma)B_K(x,z)}}
\,d\mu_x^{s_i}(\xi)\\
&=&
\lim_{j\to\infty}A_j(z)\lim_{s_i\to\delta(\Gamma)^+}
\frac{1}{c(s_i)}
\sum_{\scriptsize\begin{array}{c}g\in\Gamma\\gO\in U_j\end{array}}
\frac{\f(D(gO)) e^{-s_id(x,gO)}}{\mu_O(U_j(z))
e^{-\delta(\Gamma)B_K(x,z)}}.
\end{eqnarray}
Moreover,
\begin{eqnarray*}
& &
\frac{\f(D(gO))e^{-s_id(x,gO)}}
{\mu_O(U_j(z))e^{-\delta(\Gamma)B_K(x,z)}}\\
&=&
\frac{\f(D(gO))}{\mu_O(U_j(z))}\cdot
\frac{e^{-s_i(d(x,gO)-d(O,gO))}}{e^{-\delta(\Gamma)B_K(x,z)}}
\cdot e^{-s_id(O,gO)}\\
&=&
\frac{\f(D(gO))}{\mu_O(U_j(z))}\cdot
\frac{e^{-\delta(\Gamma)(d(x,gO)-d(O,gO))}}{e^{-\delta(\Gamma)B_K(x,z)}}\cdot
\frac{e^{-s_i(d(x,gO)-d(O,gO))}}{e^{-\delta(\Gamma)(d(x,gO)-d(O,gO))}}
\cdot e^{-s_id(O,gO)}.
\end{eqnarray*}

From the definition of the Busemann function $B_K$, it follows that
for all
$z\in\partial\M H^k$ and $p\in\M H^k$
\begin{equation}\label{eq3_7}
\lim_{\xi\to z}\left(d(p,\xi)-d(O,\xi)\right)=B_K(p,z).
\end{equation}
Therefore, there exist two sequences $\{E_j^+\}$ and $\{E_j^-\}$,
converging to $1$ as $j\to\infty$, and such that for all $g\in\Gamma$
with $g O\in U_j(z)$
$$
E_j^-\leq
\frac{e^{-\delta(\Gamma)(d(x,gO)-d(O,gO))}}{e^{-\delta(\Gamma)B_K(x,z)}}
\leq
E_j^+.$$
Moreover, equation $(\ref{eq3_7})$ implies that the term
$(d(x,gO)-d(O,gO))$ is bounded, so for each $j$
$$
\lim_{s_i\to\delta(\Gamma)^+}
\frac{e^{-s_i(d(x,gO)-d(O,gO))}}{e^{-\delta(\Gamma)(d(x,gO)-d(O,gO))}}
=1$$
uniformly on $U_j$. Whence, since $\f$ is positive,
\begin{eqnarray}
& &
\lim_{j\to\infty}A_j(z)\lim_{s_i\to\delta(\Gamma)^+}
\frac{1}{c(s_i)}
\sum_{\scriptsize\begin{array}{c}g\in\Gamma\\gO\in U_j\end{array}}
\frac{\f(D(gO)) e^{-s_id(x,gO)}}{\mu_O(U_j(z))
e^{-\delta(\Gamma)B_K(x,z)}}\\
&\leq&
\lim_{j\to\infty}A_j(z)E_j^+\lim_{s_i\to\delta(\Gamma)^+}
\frac{1}{c(s_i)}
\sum_{\scriptsize\begin{array}{c}g\in\Gamma\\gO\in U_j\end{array}}
\frac{\f (D(gO))}{\mu_O(U_j(z))}\cdot e^{-s_id(O,gO)}\\
&=&
\lim_{j\to\infty}A_j(z)E_j^+\lim_{s_i\to\delta(\Gamma)^+}
\int_{U_j(z)}\frac{\f(D(x))}{\mu_O(U_j(z))}
\,d\mu_O^{s_i}(x)
\end{eqnarray}
and similarly for $E_j^-$. Since for $\mu_O$-almost all
$z$ we have $\lim_j A_j(z)=1$, and since $E_j^\pm\to1$,
we get that $\mu_O$-almost everywhere
\begin{eqnarray}
& &\lim_{j\to\infty}A_j(z)E_j^\pm\lim_{s_i\to\delta(\Gamma)^+}
\int_{U_j(z)}\frac{\f(D(\xi))}{\mu_O(U_j(z))}
\,d\mu_O^{s_i}(\xi)\\
&=&\lim_{j\to\infty}\lim_{s_i\to\delta(\Gamma)^+}
\int_{U_j(z)}\f(D(\xi))\frac{1}{\mu_O(U_j(z))}\,d\mu_O^{s_i}(\xi)\\
&=&\lim_{j\to\infty}\lim_{s_i\to\delta(\Gamma)^+}
\intK\psi_j^{z,\eta_O}(\xi)\f(D(\xi))\,d\mu_O^{s_i}(\xi)\\
&=&\lim_{j\to\infty}\lim_{s_i\to\delta(\Gamma)^+}
\intKN\psi_j^{z,\eta_O}(\xi)\f(y)\,d\eta_O^{s_i}(\xi,y)\\
&=&\lim_{j\to\infty}\intKN\psi_j^{z,\eta_O}(\xi)\f(y)\,d\eta_O(\xi,y)\\
\label{eq_9}&=&\intN\f(y)\,d\alpha_z(y).
\end{eqnarray}
Finally, from equations~$(\ref{eq_6})-(\ref{eq_9})$, $\mu_O$-almost
 everywhere
 we get
$$
\intN\f(y)\,d\alpha_z(y)\leq
\intN\f(y)\,d\wt\alpha_z(y)\leq
\intN\f(y)\,d\alpha_z(y)
$$
and the claim follows. \qed

In particular, Theorem~\ref{t_adisint} implies that
 the measures $\eta_x$'s are $\rho$-equivariant. Indeed,
since the measures $\{\eta_x^s\}$ are $\rho$-equivariant, and since
the push-forward is continuous for the weak convergence, for all $x\in
\M H^k$ and $\gamma\in\Gamma$
$$\eta_{\gamma
  x}^{s_i}=(\gamma,\rho(\gamma))_*\eta_x^{s_i}
\rightharpoonup(\gamma,\rho(\gamma))_*\eta_x.$$
Moreover, since $\pi_{K*}(\eta_x)=\mu_x$, each measure $\eta_x$
is concentrated in $\partial \M H^k\times
\wb{\M H}^n$.

 \begin{lemma}\label{l4_4.3}
   The family $\{\alpha_z\}$ is $\rho$-equivariant, that is, for all
   $\gamma\in\Gamma$ and $\mu_O$-almost all $z\in\partial{\M H}^k$
$$\alpha_{\gamma z}=\rho(\gamma)_*\alpha_z.$$
 \end{lemma}
\proof
From point~$(\ref{psm_3})$ of
Theorems~\ref{t_psm} and the
the $\rho$-equivariance of the $\eta_x$'s, it follows that
for all $\f\in C(\hKN)$
\begin{eqnarray*}
&&  \intbKN\f(\gamma z,y)\,d\alpha_{\gamma z}(y)\,d\mu_O(z)\\
&=&\intbKN\f(z,y)\,d\alpha_z(y)\,d\mu_{\gamma O}(z)\\
&=&\intKN\f(z,y)\,d\eta_{\gamma O}=\intKN\f(\gamma z,\rho(\gamma)y)\,d\eta_O\\
&=&\intbKN\f(\gamma z,\rho(\gamma) y)\,d\alpha_z(y)\,d\mu_O(z)\\
&=&\intKN\f(\gamma z,y)\,d(\rho(\gamma)_*\alpha_z)(y)\,d\mu_O(z).
\end{eqnarray*}
Whence, the measures $\alpha_{\gamma z}$ and $\rho(\gamma)_*\alpha_z$
equal $\mu_O$-almost everywhere.
\qed

Now, for each $y\in\wb{\M H}^n$ let $\nu_y$ be the visual measure on
$\partial \M H^n$ centred at $y$. More precisely, choose the disc
model of $\M H^n$ whose centre is $O$ and  let $\nu_O$ be the
standard probability measure on $S^{n-1}\simeq \partial\M H^n$.
Then, for all $g\in\isomn$ define
$$\nu_{gO}=g_*\nu_O.$$
This definition is not ambiguous because $\nu_O$ is ${\rm
  Stab}(O)$-invariant, where ${\rm Stab}(O)=\{g\in\isomn:\, g(O)=O\}$.
For $y\in\partial\M H^n$ simply define $\nu_y=\delta_y$.

For all $z\in\partial \M H^k$ define $\lambda_z$ a measure on
$\partial \M H^n$ by $\lambda_z=\alpha_z*\{\nu_y\}$.
That is, for all $\f\in C(\partial \M H^n)$
$$\intbN\f\,d\lambda_z=
\intN\left(\intbN\f(\theta)\,d\nu_y(\theta)\right)\,d\alpha_z(y).$$
Note that such an integral is well-defined because, for any \f, the function
$y\mapsto\intbN\f\,d\nu_y$ is continuous in $y$.
\begin{remark}
  Since the measure $\nu_y$ depends continuously on $y$, by
  Lemma~\ref{l2.0_2.9}, the
  convolution with the family of visual measures is weakly continuous.
\end{remark}
We show now that $\{\lambda_z\}_{z\in\partial\M H^k}$ is a
family of developing measures for $\rho$.
\begin{lemma}
  The family $\{\lambda_z\}$ is $\rho$-equivariant, that is, for
  $\mu_O$-almost all $z\in\partial\M H^k$ and all
  $\gamma\in\Gamma$, we have $\lambda_{\gamma z}=\rho(\gamma)_*\lambda_z.$
\end{lemma}
\proof
By Lemma~\ref{l4_4.3}, for $\mu_O$-almost all $z\in\partial\M H^k$ and
all $\f\in C(\partial \M H^n)$
\begin{eqnarray*}
 \intbN\!\!\f\,d\lambda_{\gamma z}\!\!\!
&=\!\!\!\!&\intN\!\intbN\!\!\f(\theta)\,d\nu_y(\theta)\,d\alpha_{\gamma z}(y)
=\!\intN\!\intbN\!\!\f(\theta)\,d\nu_{\rho(\gamma)y}(\theta)\,d\alpha_z(y)\\
&=\!\!\!\!&\intN\intbN\f(\rho(\gamma)\theta)\,d\nu_y(\theta)\,d\alpha_z(y)\\
&=\!\!\!\!&\intbN\f\circ\rho(\gamma)\,d\lambda_z=\intbN\f\,d(\rho(\gamma)_*\lambda_z).
\end{eqnarray*}
\qed
\begin{lemma}
  For any $\f\in C(\partial\M H^n)$, the function
$z\mapsto \intbN \f(\theta)\,d\lambda_z(\theta)$ is $\mu_O$-integrable.
\end{lemma}
\proof
This follows directly from the definition of $\lambda_z$, because the
family $\{\alpha_z\}$ disintegrates $\eta_O$, and $\mu_O=\pi_{K*}\eta_O$.
\qed
\begin{lemma}
For $\mu_O$-almost all $z\in\partial\M H^k$, we have $||\lambda_z||=1$.
\end{lemma}
\proof
For all $x\in\M H^k$ and $y\in\partial \M H^n$, the measures $\eta_x$
and $\nu_y$ are positive. Then the measures $\alpha_z$'s are positive,
and this implies that the measures $\lambda_z$'s are positive.
Thus
\begin{eqnarray*}
  ||\lambda_z||=\intbN1\,d\lambda_z=
\intN\intbN1\,d\nu_y(\theta)\,d\alpha_z(y)
=\intN1\,d\alpha_z=||\alpha_z||=1
\end{eqnarray*}
$\mu_O$-almost everywhere because $\{\alpha_z\}$ disintegrates $\eta_O$.
\qed

Therefore, the family $\{\lambda_z\}$ satisfies properties
$(\ref{devm_1})-(\ref{devm_4})$ of Definition~\ref{d_devm}. So it
  is a family of developing measure for $\rho$, and the proof of
  Theorem~\ref{t_3} is complete.\qed


\section{Sequence of \e-natural maps}\label{s2.1_5}

In this section we prove Theorem~\ref{p_5}, showing
 that the natural map constructed in Sections~\ref{s_3}
and~\ref{s_4} is the limit of a sequence of \e-natural maps.
We keep the notation of Sections~\ref{s_2}-\ref{s_4}.
Through this section we
suppose that $\Gamma$ diverges at $\delta(\Gamma)$; by
Lemma~\ref{l_nic}, this is the case
if $\M H^k/\Gamma$ has finite volume.

Let $A\subset \M H^k$ be the Dirichlet domain of $O$. The set $A$ is a
fundamental domain for $\Gamma$ containing $O$.
Let $\sigma$ be any Borel probability measure on $A$.

\begin{defi}
  For each $s>\delta(\Gamma)$ and $x\in \M H^k$ we define $m_x^s$ a
  positive Borel measure on $\M H^k$ by
$$m_x^s=\frac{1}{c(s)}\sum_{\gamma\in\Gamma}e^{-sd(x,\gamma O)}\gamma_*\sigma.$$
\end{defi}

Note that if $\sigma=\delta_O$, then $m_x^s=\mu_x^s$.

\begin{lemma}\label{l_2.1_5.3}
Let $\{\gamma_n\}$ be any sequence of elements of $\Gamma$.
If $\gamma_n(O)$ converges to a point $\theta\in\partial\M H^k$, then
for any sequence $\{x_n\}\in A$, the sequence $\{\gamma_n(x_n)\}$ converges to
$\theta$.
\end{lemma}
\proof Suppose the contrary. Since $\wb{\M H}^k$ is compact, up to
passing to a sub-sequence, we can suppose that $\{\gamma_n(x_n)\}$
converges to a point $\zeta\neq\theta$ in $\wb{\M H}^k$. Then, the
geodesics joining $\gamma_n(x_n)$ to $\gamma_n(O)$ accumulate near the geodesic
between $\zeta$ and $\theta$. This can not happen because $\Gamma$ is
discrete.\qed

\begin{teo}
For each $x\in\M H^k$, if $\mu_x$ denotes the Patterson-Sullivan
measure constructed as in Theorem~\ref{t_psm}, then
$$m_x^s\rightharpoonup\mu_x$$
in $\wb{\M H}^k$, when $s\to\delta(\Gamma)^+$.
\end{teo}
\proof We have to show that for each $\f\in C(\wb{\M H}^k)$,
$\int\f\,dm_x^s\to\int\f\,d\mu_x$. Let $\f\in C(\wb{\M H}^k)$. Since
$\mu_x^s\to \mu_x$, we will have finished by proving that
$$\lim_{s\to\delta(\Gamma)^+}\Big|\intK\f\,d m_x^s-\intK\f\,\mu_x^s\Big|=0.$$

Let $C>0$ be a small constant and let $A_1\subset A$ be a compact set such
that $O\in A_1$ and $\sigma(A\setminus A_1)<C$. Since the supports
of the measures $m_x^s$ and $\mu_x^s$ are contained in $\M H^k$, and
since $O\in A$, we have
\begin{equation}\label{eq3.0_boh}
\intK\f\,d(m_x^s-\mu_x^s)=
\sum_{\gamma\in\Gamma}
\int_{\gamma(A_1)}\f\,d(m_x^s-\mu_x^s)+
\sum_{\gamma\in\Gamma}
\int_{\gamma (A\setminus A_1)}\f\,dm_x^s.
\end{equation}
Looking at the second summand,
\begin{eqnarray}
& & \Big| \sum_{\gamma\in\Gamma}
\int_{\gamma (A\setminus A_1)}\f\,dm_x^s\Big|\\
&\leq&
\max(|\f|) \frac{1}{c(s)}\sum_{\gamma\in\Gamma}
\int_{\gamma (A\setminus A_1)}e^{-sd(x,\gamma O)}\,d\gamma_*\sigma(\xi)\\
\label{e2.2_18}
&=&
\max(|\f|) \frac{1}{c(s)}\sum_{\gamma\in\Gamma}
e^{-sd(x,\gamma O)}\sigma(A\setminus A_1)\leq
C\cdot\max(\f)\cdot||\mu_x^s||.
\end{eqnarray}

We estimate now the first summand.
Since $c(s)\to\infty$, for any finite subset $\Gamma_1$ of $\Gamma$
\begin{equation}\label{e2.1_17}
\lim_{s\to \delta(\Gamma)^+}
\Big|\sum_{\gamma\in\Gamma_1}\int_{\gamma(A_1)}\f\,d m_x^s\Big|
=
\lim_{s\to \delta(\Gamma)^+}
\Big|\sum_{\gamma\in\Gamma_1}\int_{\gamma(A_1)}\f\,d \mu_x^s\Big|
=0.
\end{equation}

Moreover, since $\f$ is continuous, it is uniformly continuous on
$\wb{\M H}^k$ for any metric that induces the usual topology on
$\wb{\M H}^k$ (recall that the hyperbolic metric of $\M H^k$ is not a metric on
$\wb{\M H}^k$).
Therefore, Lemma~\ref{l_2.1_5.3} implies that, except for a finite
number of elements of $\Gamma$, we have
\begin{eqnarray}\label{e2.1_18}
  |\f(\gamma O)-\f(\gamma\xi)|<C
\end{eqnarray}
independently on $\gamma$ and on $\xi\in A_1$.
Let $\Gamma_1$ be a finite subset of $\Gamma$ such that
$(\ref{e2.1_18})$ holds for $\gamma\in\Gamma\setminus\Gamma_1$.
Then,
\begin{eqnarray*}
& &
\Big|\sum_{\gamma\in\Gamma\setminus\Gamma_1}
\int_{\gamma (A_1)}\f\,d(m_x^s-\mu_x^s)\Big|\\
&=&
\Big|\frac{1}{c(s)}\sum_{\gamma\in\Gamma\setminus\Gamma_1}
e^{-sd(x,\gamma O)}\Big(
\int_{\gamma(A_1)}\f(\xi)\,d\gamma_*\sigma(\xi)\\ & &
\hskip30ex -\, \f(\gamma O)\big(\sigma(A_1)+\sigma(A\setminus A_1)\big)
\Big)\Big|\\
&\leq& C\intK|\f|\,d\mu_x^s+
\frac{1}{c(s)}\sum_{\gamma\in\Gamma\setminus\Gamma_1}
e^{-sd(x,\gamma O)}\int_{A_1}
|\f(\gamma\xi)-\f(\gamma O)|
\,d\sigma(\xi)
\\
&\leq& C\intK|\f|\,d\mu_x^s+
\frac{C}{c(s)}\sum_{\gamma\in\Gamma\setminus\Gamma_1}
e^{-sd(x,\gamma O)}\cdot\sigma(A_1)
\leq C\intK(1+|\f|)\,d\mu_x^s.
\end{eqnarray*}

Whence the claim follows, combining with
$(\ref{eq3.0_boh}),\ (\ref{e2.2_18})$ and $(\ref{e2.1_17})$,
since $C$ can be chosen arbitrarily small.\qed

Now we proceed as in Sections~\ref{s_3} and~\ref{s_4}.
Namely, we fix a $\rho$-equivariant map $D$, we define measures
$n_x^s=m_x^s\times\{\delta_{D(z)}\}$, and we chose
a sequence $s_i\to\delta(\Gamma)^+$ such that $n_O^{s_i}$ converges to a
measure $n_O$. We disintegrate $n_O$ as $n_O=\mu_O\times\{a_z\}$, and
we define $n_x=\mu_x\times\{a_z\}$. As in Theorem~\ref{t_adisint} one
can show that $n_x^{s_i}\rightharpoonup n_x$. We define then
$b_x^s=(D_*m_x^s)*\{\nu_y\}=\pi_{N*}(n_x^s)*\{\nu_f\}$, and
$b_x=\pi_{N*}(n_x)*\{\nu_f\}$.

Finally, for each $s>\delta(\Gamma)$ we set $s=(1+\e)\delta(\Gamma)$
and  we define maps $F(x)=\B(b_x)$ and
$F_\e(x)=\B(b_x^s)$. The map $F$ is a natural
map, in the sense that it has the properties $(1)-(4)$ of Theorem~\ref{t_1}.
The maps $F_\e$ the \e-natural maps constructed in~\cite{bcg99}, so
they are smooth, $\rho$-equivariant, and
for all $p\geq 3$ and $\e>0$,\  $\displaystyle{
{\rm  Jac}_pF_\e(x)\leq\left(\frac{(1+\e)\delta(\Gamma)}{p-1}\right)^p}.$

  \begin{prop}\label{p_3.0_6.4}
    The maps $F_{\e_i}$ punctually converge to the map $F$.
  \end{prop}

\proof
From the weak continuity of the
push-forward and from Lemma~\ref{l2.0_2.9}, we get
$$b_x^{s_i}\rightharpoonup b_x.$$
Then, the claim follows from Point~\ref{eq_bar3} of Lemma~\ref{l_bar}.\qed

\begin{remark}\label{r_6.0_R}
The weak convergence of $b_x^s$ to $b_x$ is enough to
  prove stronger convergences. For example, it can be shown that the
  derivatives of $F_{\e_i}$ converges the ones of $F$, whence one gets
  that the convergence of the \e-natural maps is locally uniform
  (see~\cite{frkl06} for details).
\end{remark}

\section{Rigidity of representations}\label{s_R}
In this section we give a proof of Theorem~\ref{t_6.0_rigidity},
referring the reader to~\cite{frkl06} for a fully
detailed discussion on the matter.

\ \\
\noindent{\em Proof of Theorem~\ref{t_6.0_rigidity}.}
Let $M$ be a complete hyperbolic $k$-manifold of finite volume and let
$\rho:\pi_1(M)\to \isomn$ be a representation. We consider $\M H^k$ as the
universal cover of $M$, and we identify $\pi_1(M)<\isomk$ with the group
of deck transformations of $\M H^k\to M$.

We denote by a {\em pseudo-developing map} any piecewise smooth,
$\rho$-equivariant map $D:\M H^k\to\M H^n$. The volume of a pseudo
developing map $D$ is defined by integrating the pull-back of the volume
form of $\M H^n$ on a fundamental domain for $M$. Equivalently, if
$\omega_N$ denotes the volume form of $\M H^n$, by equivariance, the
form $D^*\omega_N$ descends to $M$ and we set
$$\vol(D)=\int_M D^*\omega_N.$$

Now, let us suppose that $M$ is compact. In this case we define the
volume of $\rho$ as the infimum of the volumes of all the
pseudo-developing maps for $\rho$:
$$\vol(\rho)=\inf_D\vol(D).$$

Note that the volume of any elementary representation vanishes, so
that we can suppose $\rho$ to be non-elementary.

The existence of a natural map immediately gives
the inequality; indeed  Theorem~\ref{t_1} point~\ref{t1_p3},
together with the fact that
$\delta(\pi_1(M))=k-1$, tells us that the volume of a natural map is
less than $\vol(M)$.
Moreover, if $\rho$ has maximum volume, then from Theorem~\ref{t_1}
point~\ref{t1_p4} it follows that in each point the differential of a
natural map is an isometry. Since $\vol(\rho)$ is maximal,
we deduce that the image of a natural map is a locally minimal sub-manifold
of $\M H^n$. The claim now follows because a locally isometric and
locally minimal immersion from a Riemannian manifold to another is
totally geodesic.

\ \\

In the non-compact case, some problem arises. From now on we
suppose that $M$ is non-compact. If we keep the above definition of
volume of $\rho$ we get that the volume of any representation
vanishes. This is because any non-compact manifold collapses to a
spine, which is a codimension-one object. To avoid such a pathology we
need to require that a pseudo-developing map has a nice behaviour on
the cusps.
For $G$ a group of isometries, we denote by Fix$(G)$ the set of fixed
points of $G$ (including the points at infinity). If $C$ is a cusp of
$M$, it is readily checked that $\pi_1(C)$ is Abelian and that
Fix$(\pi_1(C))$ consists of a unique point of $\partial \M H^k$.

We say that a pseudo-developing map $D$ {\em properly ends} if for any
cusp $C$ of $M$, if
$\xi=\textrm{Fix}(\pi_1(C))$ and $\alpha(t)$ is a
geodesic ray ending at   $\xi$, then all limit points of
$D(\alpha(t))$ lie either in
$\textrm{Fix}(\rho(\pi_1(C)))\subset\wb{\M H}^n$
 or in a finite union
  of $\rho(\pi_1(C))$-invariant geodesics. It easy to see that
  properly ending pseudo-developing maps always exist.

Now we can define the volume of $\rho$ by taking the infimum of the
volumes of properly ending pseudo-developing maps.
$$\vol(\rho)=\inf_{D\textrm{ properly ending}}\vol(D).$$

We would like to use a natural map as in the compact case, but
the problem now is that a natural map given by Theorem~\ref{t_1} in
general does not end properly.

Nevertheless, we can use the \e-natural maps.
Indeed, it can be shown
(see~\cite{frkl06}) that any \e-natural map
constructed as in Section~\ref{s2.1_5} properly ends, and this gives
immediately the inequality.
Moreover, by Proposition~\ref{p_3.0_6.4} and Remark~\ref{r_6.0_R} we
have that the volumes of the \e-natural maps converge to the volume
of a natural map. Thus, if the volume of $\rho$ is maximal, then the
volume of a natural map is maximal and we conclude as in the compact
case.\qed


\section{Measurable extension of natural maps }\label{s_5}
This section is devoted to proving Theorem~\ref{t_4}.
We keep here the notation of previous sections. In particular we
recall that $\{\mu_x\}_{x\in\M H^k}$ is the family of
Patterson-Sullivan measures, and that
$\{\lambda_z\}_{z\in\partial \M H^k}$ is a family of developing
measures.

\begin{defi}
Let $\{\gamma_i O\}$ be a sequence in the $\Gamma$-orbit of $O$. We
say that $\gamma_i O$ {\em conically} converges to $\omega\in\partial
\M H^k$ if $\gamma_iO\to \omega$ and there exists a geodesic $\sigma$,
ending at $\omega$, such that the distance of $\gamma_iO$ from
$\sigma$ is bounded.
The {\em conical limit set} of $\Gamma$, denoted by $\Lambda_c$,
is the set of the limits of conically converging sequences in the
$\Gamma$-orbit of $O$.
\end{defi}

Clearly, the conical limit set is a sub-set of the limit set of
$\Gamma$.
In order to prove that the natural maps extend to the boundary we need
the following result.
\begin{teo}\label{t_4.0_6.2}
  For each $f\in L^1(\partial \M H^k,\mu_O)$ there exists a
   set $Z$ with $\mu_O(Z)=0$ such that for all
   $\omega\in\Lambda_c\setminus Z$, and for any sequence
   $\{\gamma_i\}\subset \Gamma$ such that $\gamma_iO$ conically
   converges to $\omega$
$$\intbK f(\theta)\,d\mu_{\gamma_iO}(\theta)\to f(\omega).$$
\end{teo}

Before proving Theorem~\ref{t_4.0_6.2} we show how it implies
Theorem~\ref{t_4}.

\proof[Proof of Theorem~\ref{t_4}].
First, we prove the existence part.
  By definition (see Section~\ref{s_2}) we
have
\begin{eqnarray*}
&  & \call{B}_{\beta_x}(y)=\intbN B_N(y,\theta)\,d\beta_x(\theta)=
\intbK\intbN B_N(y,\theta)\,d\lambda_z(\theta)\,d\mu_x(z)\\
&=& \intbK\call B_{\lambda_z}(y)\,d\mu_x(z)
\end{eqnarray*}
and a similar formula holds for the derivatives of $\call
B_{\beta_x}$.
Therefore, for each $y$ there exists a
$\mu_O$-negligible set $Z\subset\partial \M H^k$
such that for all $\omega\in\partial \M H^k\setminus Z$ we have
$$\lim_{x\to\omega}\call B_{\beta_{x}}(y)=B_{\lambda_\omega}(y)$$
Where ``$\lim_{x\to \omega}$'' means ``for any sequence $\{x_i\}$ in
the $\Gamma$-orbit of $O$, conically converging to $\omega$...''. The
same statement holds for the derivatives of $\call
B_{\beta_{x}}$.
Now, let $Y$ be a countable dense subset of $\M H^n$. Then (since
countable unions of negligible sets are negligible,) there exists a
$\mu_O$-negligible set $W\subset\partial \M H^k$ such that the above
limit holds for all $\omega\in\partial\M H^k\setminus W$, all $y\in Y$
and all the derivatives of $\call B_{\beta_x}$.
It follows that the barycentre of $\beta_x$, that is the unique point
of minimum of $\call B_{\beta_x}$, converges to the barycentre of
$\lambda_{\omega}$, which is well-defined because
$\lambda_\omega$ is not the sum of two Dirac deltas with equal weights.
Therefore, if $F$ denotes the natural map constructed using the family
$\{\lambda_z\}$, setting $\wb F(\omega)=\B(\lambda_\omega)$, we have that
for $\mu_O$-almost all $\omega$, for any sequence $\{\gamma_i O\}$
conically converging to $\omega$

\begin{equation}\label{e_7.0_21}
\lim_{\gamma_i O\to\omega}F(\gamma_i O)=\wb F(\omega)=\B(\lambda_\omega).
\end{equation}

The map $\wb F$ is measurable because it can be viewed as a limit of
continuous functions. Finally, it is readily checked that
$F(\gamma_i O)$ and $F(\gamma_i x)$ have the same limit, and
 this completes the proof of the existence part.
It remains now to prove the last part of Theorem~\ref{t_4}.

Given the maps $\wb F_1$ and $\wb F_2$, we construct the corresponding
natural maps $F_1, F_2$. For $i=1,2$ and for
almost all $\omega$ in the conical limit set of $\Gamma$, if
$\{\gamma_nO\}$ is a sequence conically converging to $\omega$, then
by~(\ref{e_7.0_21}),
$F_i(\gamma_nO)\to\B(\wb F_i(\omega))$. By equivariance we have
$$F_i(\gamma_nO)=\rho(\gamma_n)F_i(O)\qquad i=1,2.$$

Up to pass to a sub-sequence, either $\rho(\gamma_n)$ converges to an
isometry $\psi$ of $\M H^n$, or the limit of $\rho(\gamma_n)y$ belongs
to $\partial\M H^n$ and does not
depend on $y\in\M H^n$
(see for example~\cite{Kap:libro}).
In the former case, we get $\psi(F_i(O))=\B(\wb
F_i(\omega))$; but also, for all $\gamma\in\Gamma$ we have
$$\psi(\rho(\gamma)F_i(O))=\psi(F_i(\gamma
O))=\lim\rho(\gamma_n)F_i(\gamma O).$$
The sequence $\gamma_n\gamma O$ conically converges to $\omega$, since
$\gamma_n O$ does. Then, by equivariance and by(~\ref{e_7.0_21}), we get
$$\lim\rho(\gamma_n)F_i(\gamma O)=\lim F_i(\gamma_n(\gamma O))=\B(\wb
F_i(\omega))$$
whence $\psi(\rho(\gamma)F_i(O))=\psi(F_i(O))$ and thus
$\rho(\gamma)F_i(O)=F_i(O)$, contradicting the fact that the
image of $\rho$ is non-elementary.
Therefore, we are in the latter case and in particular
$$\B(\wb F_1(\omega))=\lim \rho(\gamma_n)F_i(O)=\lim
\rho(\gamma_n)F_2(O)=\B(\wb F_2(\omega))\in\partial\M H^n$$
but this is possible if and only if
$\wb F_i(\omega)=\delta_{\B(F_i(\omega))}$. So $\wb F_1$ and $\wb F_2$
are ordinary functions with values in $\partial \M H^n$,
and they coincide almost everywhere.
\qed

\proof[Proof of Theorem~\ref{t_4.0_6.2}]. The map
$x\mapsto \intbK f(\theta)\,d\mu_x(\theta)$ can be viewed as the
harmonic extension of $f$, and one can
prove its convergence to $f$ along cones and almost everywhere
by using standard techniques of harmonic analysis. We give a
proof for completeness.

For the whole proof, we work in the half space model $\M H^k=\M
R^{k-1}\times\M R^+$, using the following notation. For a point
$x\in\M H^k$, we denote by $(x',x_k)\in \M R^{k-1}\times\M R^+$
its coordinates in the half-space model, by the symbol $|x|$
we denote the Euclidean norm of $x$ in the model, and $B(\omega,r)$
will denote the ball of centre $\omega$ and Euclidean radius $r$.
In the half-space
model, setting $\delta=\delta(\Gamma)$ we have
$$e^{-\delta B_K(x,\xi)}=
\left(\frac{x_k(1+|\xi|^2)}{|\xi-x'|^2+x_k^2}\right)^\delta$$
and for all $\gamma\in\Gamma$, if $x=\gamma O$

$$1=||\mu_O||=||\gamma_*\mu_O||=||\mu_{\gamma O}||=
\intbK \left(\frac{x_k(1+|\xi|^2)}{|\xi-x'|^2+x_k^2}\right)^\delta.$$

We will work on a fixed ball of centre $0$ and radius $R$ of $\M
R^{k-1}$. This is not restrictive because proving convergence
almost everywhere for all balls is equivalent to prove convergence
almost everywhere.

For all $\omega\in\M R^{k-1}$, we denote by $C_\omega(\alpha)$ the
vertical cone in $\M H^{k-1}\times\M R^+$ of vertex $\omega$ and
emi-angle $\alpha$. For any non-negative $g\in L^1(\partial \M H^k,\mu_O)$,
$\omega\in\partial \M H^k$, $\alpha\in(0,\pi/2)$ we define the
maximal operator $M_\alpha g(\omega)$ by
$$M_\alpha g(\omega)=\sup_{\gamma\in\Gamma,\
\gamma O\in C_\omega(\alpha)}\intbK
g(\xi)\,d\mu_{\gamma O}(\xi)$$
and the so-called Hardy Littlewood operator $N g(\omega)$ by
$$Ng(\omega)=\sup_{r>0}\frac{1}{\mu_O(B(\omega,r))}\int_{B(\omega,r)}g(\xi)\,d\mu_O(\xi).$$

From now on, the symbol $c$ will denote a generic constant, and
different occurrences may denote different constants. If not specified,
the constants do not depend on the other quantities we are considering.

\begin{lemma}\label{l_4.0_6.3}
  There exists a constant $c$ such that for every point $\omega$  of the limit
  set of $\Gamma$ and any $r>0$
$$\mu_O(B(\omega,r))\leq cr^\delta.$$
\end{lemma}

\proof
For all $x\in\M H^k$
$$||\mu_x||=
\intbK\left(\frac{x_k(1+|\xi|^2)}{|\xi-x'|^2+x_k^2}\right)^\delta
\geq
\intbK\left(\frac{x_k}{|\xi-x'|^2+x_k^2}\right)^\delta.
$$
Suppose now that $x$ is of the form $x=\gamma O$, with
\begin{equation}\label{eq4.0_21}
c_1r\leq |x-\omega|\leq c_2r.
\end{equation}
Then $|\xi-x'|\leq |\xi-\omega|+|x'-\omega|\leq(1+c_2)r$
for any $\xi\in B(\omega,r)$. Then
$$1=||\mu_x||=
\intbK\left(\frac{x_k}{|\xi-x'|^2+x_k^2}\right)^\delta
\geq
c\frac{1}{r^\delta}\mu_O(B(\omega,r))
$$
and the claim holds for such points. Let $a>1$ and for
$j\in\M Z$, consider the set $A_j=\{x\in\wb{\M H}^k :\
|x-\omega|\in [a^j,a^{j+1})\}$. If for a certain $j\in\M Z$ the set
$A_j$ contains a point of the
$\Gamma$-orbit of $O$, then for $r\in[a^j,a^{j+1}]$
 $$\frac{r}{a}\leq\frac{a^{j+1}}{a}=
a^j\leq x\leq a^{j+1}=a^{j}a\leq r a$$
and inequalities~$(\ref{eq4.0_21})$ hold with $c_2=a=1/c_1$.
Let now $r>0$. Since $\omega$ lies on the limit set, there exists
$j$ with $a^j\leq r$ and such that $A_{j}$ intersects the $\Gamma$-orbit
of $O$. Let $j_0$ be the maximum of such $j$'s. If $a^{j_0+1}\geq r$,
we have finished. Otherwise, since $\mu_O$ is concentrated on the
limit set
$$\mu_0(B(\omega,r))=\mu(B(\omega,a^{j_0+1}))\leq c
(a^{j_0+1})^\delta\leq c r^\delta.$$
This completes the proof of Lemma~\ref{l_4.0_6.3}.\qed

\begin{lemma}\label{l_4.0_6.4}
There exists a constant $c$, depending only on $\alpha$, such that for
any non-negative $g\in L^1(\partial \M H^k,\mu_O)$,
$\omega\in\partial \M H^k$, $\alpha\in(0,\pi/2)$, we have
$$M_\alpha g(\omega)\leq c N g (\omega).$$
\end{lemma}
\proof
Let $x=\gamma O$ with $\gamma\in\Gamma$ and $\gamma O\in
C_\omega(\alpha)$. As noticed above, it is not restrictive to work in
the ball $B(0,R)$.
\begin{eqnarray*}
&~&
\intbK g(\xi)\,d\mu_x(\xi)\leq (1+R^2)^\delta\intbK
g(\xi)\left(\frac{x_k}{|\xi-x'|^2+x_k^2}\right)^\delta\,d\mu_O(\xi)
\\
&=&
c\Bigg[\int_{B(\omega,x_k)}
g(\xi)\left(\frac{x_k}{|\xi-x'|^2+x_k^2}\right)^\delta\,d\mu_O(\xi)+
\\
& &+\sum_{j\geq 0}
\int_{B(\omega,2^{j+1}x_k)\setminus B(\omega,2^jx_k)}
g(\xi)\left(\frac{x_k}{|\xi-x'|^2+x_k^2}\right)^\delta\,d\mu_O(\xi)
\Bigg]\\
&\leq&
c\Bigg[\frac{1}{x_k^\delta}\int_{B(\omega,x_k)} g(\xi)\,d\mu_O(\xi)+\\
& & \hskip10ex +
\sum_{j\geq 0} \frac{1}{(c2^{2j}x_k^2)^\delta}
\int_{B(\omega,2^{j+1}x_k)\setminus B(\omega,2^jx_k)}
g(\xi)\,d\mu_O(\xi)\Bigg]\\
&\leq&
c\Bigg[\frac{1}{\mu_O(B(\omega,x_k))}\int_{B(\omega,x_k)}g(\xi)\,d\mu_O(\xi)+
\\ & & \hskip10ex + \sum_{j\geq 0}
\frac{2^{-j}}{\mu_O(B(\omega,2^{j+1}x_k))}\int_{B(\omega,2^{j+1}x_k)}
g(\xi)\,d\mu_O(\xi)
\Bigg]\\
&\leq& c N g (\omega)\left(1+\sum_{j\geq0}2^{-j}\right)\leq cNg(\omega).
\end{eqnarray*}
\qed

The constant $c$ actually depends on $\alpha$ because we used that for
$x\in C_\alpha(\omega)$ and $\xi \in B(\omega,2^{j+1}x_k)\setminus
B(\omega,2^jx_k)$ we have
 $|\xi-x'|^2+x_k^2\geq c 2^{2j}x_k^2$. It can be shown that $c$ is
 bounded by $(\tan\alpha)^{2\delta}$.

We can now finish the proof of Theorem~\ref{t_4.0_6.2}.
Since $||\mu_{\gamma O}||=1$ and $\lim_{x\to\omega}e^{-\delta B^K(x,z)}=0$
for all $z\neq \omega$, the claim is true for continuous
functions. Suppose now $f\in L^1(\partial \M H^k,\mu_O)$, and let
$f_j\to f$ be a sequence of continuous functions converging to $f$
$\mu_O$-almost  everywhere and in $L^1$. We have
\begin{eqnarray*}
\left|\intbK f(\xi)\,d\mu_{\gamma O}(\xi)-f(\omega)\right|
&\leq&
\left|\intbK f(\xi)-f_j(\xi)\,d\mu_{\gamma O}(\xi)\right|+\\&+&
\left|\intbK f_j(\xi)\,d\mu_{\gamma O}(\xi) -f_j(\omega)\right|+\\&+&
|f_j(\omega)-f(\omega)|
\end{eqnarray*}

The second summand of the second term goes to zero as $\gamma
O\to\omega$ because $f_j$ is continuous. For $\mu_O$-almost all
$\omega$ the last summand can be chosen arbitrarily small
because $f_j\to f$ $\mu_O$-almost everywhere. By
Lemma~\ref{l_4.0_6.4}, the first summand of the second term is bounded
by $cN(f_j-f)(\omega)$ on each cone $C_\alpha(\omega)$, the constant
$c$ depending on $\alpha$. The Hardy Littlewood operator is bounded
from $L^1$ to $L^{1,\infty}$ (see for example~\cite{Ste:libro70}), that is,
for all $g\in L^1$ and $\epsilon>0$
$$\mu_O\left(\{\omega\in\partial \M H^k:\ |Ng(\omega)|\geq \epsilon\}\right)
\leq\frac{c||g||_{L^1}}{\epsilon~}$$

Let $A_j^\epsilon=\{\omega\in\partial \M H^k:\ |N(f_j-f)(\omega)|\geq
\epsilon\}$. Since $f_j\to f$ in $L^1$, for all $\e$ the measure
$\mu_O(A_j^\epsilon)$ goes to zero. This is equivalent to say that the
characteristic function $\chi_{A_j^\epsilon}\to 0$ in $L^1$. Then, up
to pass to sub-sequences, $\chi_{A_j^\epsilon}\to 0$ $\mu_O$-almost
everywhere, that is, for $\mu_O$-almost all $\omega$ the quantity
$|N(f_j-f)(\omega)|\leq \epsilon$ eventually on $j$. We have so proved
that, for each $\alpha\in (0,\pi/2)$ there exist a negligible set
$Z^\epsilon_\alpha$ such that, for all $\omega\in\partial \M
H^k\setminus Z^\epsilon_\alpha$,
the quantity $|\intbK f(\xi)\,d\mu_{\gamma O}(\xi)-f(\omega)|$ is
small than or equal to $\epsilon$ as $\gamma O$ converges to $\omega$ through
$C_\alpha(\omega)$. The thesis now follows setting
$$Z=\bigcup_{
\begin{array}{c}
\alpha\in\M Q\cap(0,\pi/2)\\
0<\epsilon\in\M Q
\end{array}}
Z^\epsilon_\alpha.$$
\qed


\section{Measurable Cannon-Thurston maps}\label{s_7}
In this sections we study existence, uniqueness and convergence of
measurable Cannon-Thurston maps.

We keep here the notations of previous sections. In particular, if not
specified, $\Gamma$ is a discrete group of $\isomk$,
$\rho:\Gamma\to\isomn$ is a representation whose image is not
elementary, $\{\mu_x\}$ is the family of Patterson-Sullivan measure, and
$\{\lambda_z\}_{\z\in\partial \M H^k}$ is a family of developing
measures for $\rho$.

The following lemma collects some ergodic properties of
$\Gamma$ that we need in the sequel
(see~\cite[Theorem~A]{yue96},~\cite[Theorem~6.3.6]{nic:libro},
\cite{sul79,rob00} for the proof).
\begin{lemma}\label{l_7.1+2}  Any non-elementary discrete group $\Gamma$ acts
  ergodically on $\Lambda$ w.r.t. the measure
  $\mu_O$. Therefore, $\Lambda_c$ has either zero or full
  $\mu_O$-measure. Moreover, the following are equivalent:
  \begin{enumerate}
  \item $\Lambda_c$ has full measure.
  \item $\Gamma$  diverges at $\delta(\Gamma)$.
  \item The geodesic flow is ergodic.
  \item $\Gamma$ acts ergodically on
   $\partial \M H^k\times\partial \M H^k$ w.r.t. $\mu_O\times\mu_O$.
  \end{enumerate}

 \end{lemma}

To begin with, we prove a couple of lemmas we need in the sequel.
\begin{lemma}\label{l_5.0_7.3}
  The subset of $\partial \M H^k$ consisting of the points $z$ such
  that $\lambda_z$ is the sum of two Dirac deltas with equal weights
  is $\mu_O$-measurable.
\end{lemma}
\proof
Let $B(\theta,r)$ denotes the
ball of centre $\theta\in\partial \M H^n$ and radius $r$ in some
metric of $\partial \M H^n$, and let $Q$ be a countable, dense subset of
$\partial \M H^n$.

For any open set $B\subset \partial \M H^n$, the function $z\mapsto
\lambda_z(B)$ is $\mu_O$-measurable. It follows that the function
$$z\mapsto\inf_{0<r_1,r_2\in\M Q}\left(\sup_{\theta_1,\theta_2\in
  Q}\lambda_z\left(B(\theta_1,r_1)\cup B(\theta_2,r_2)\right)\right)$$
is $\mu_O$-measurable. The pre-image of $1$, which therefore is a
  $\mu_O$-measurable set, is the set of points $z$
  such that the support of $\lambda_z$ contains at most two points.
Similarly, the following sets are $\mu_O$-measurable
$$\{z\in\partial \M H^k:\, \lambda_z \textrm{ is concentrated on one
 point}\}$$
$$\{z\in\partial \M H^k:\, \lambda_z \textrm{ has an atom of weight }\frac{1}{2}\}$$
 and the claim follows.
\qed
\begin{lemma}
  For all
  $\mu_O$-measurable set $A\subset \partial \M H^k$ the set
$$\call O(A)=\{(x,y)\in\partial\M H^k\times\partial\M H^k:\ \exists
\gamma\in\Gamma:\, \gamma(x),\gamma(y)\in A\}$$
is $\mu_O\times\mu_O$-measurable.
\end{lemma}
\proof
Clearly, it is sufficient to show that the function
$$(x,y)\mapsto \#\{\gamma\in \Gamma:\, \gamma(x),\gamma(y)\in A\}$$
is $\mu_O\times\mu_O$-measurable.
The pre-image of $(n,\infty]$ is the set
$$\bigcup_{
  \begin{array}{c}
\gamma_1,\dots,\gamma_n\in\Gamma\\
\gamma_1\neq\dots\neq\gamma_n
\end{array}}
\Big(\big((\gamma_1(A)\times\gamma_1(A)\big)\cap\dots\cap
\big((\gamma_n(A)\times\gamma_n(A)\big)\Big)
$$
which is a countable union of countable measurable sets, and therefore
it is measurable.\qed

\begin{lemma}\label{l_5.0_7.5}
  Suppose that $\Gamma$ diverges at
  $\delta(\Gamma)$. If $\mu_O(A)>0$, then $\call O(A)$ has full measure.
\end{lemma}
\proof
Since $\mu_O(A)>0$, then $A\times A$ has positive measure. Thus,
$\call O(A)$
has positive measure because it contains $A\times A$. Moreover, it is
readily checked that $\call O(A)$ is $\Gamma$-invariant. By
Lemma~\ref{l_7.1+2} the action of $\Gamma$ on $\partial\M
H^k\times\partial\M H^k$ is ergodic, whence $\call O(A)$ has full
measure.\qed

\begin{teo}\label{t_lemmaF3}
Let $\Gamma<\isomk$ be a discrete group which diverges at
$\delta(\Gamma)$.  Let $\rho:\Gamma\to\isomn$ be a
representation whose image is not elementary and let
$\{\lambda_z\}_{z\in\partial \M H^k}$ be a family of developing measures
for $\rho$. Then,  for almost all $z$, the measure
$\lambda_z$ is not the sum of two Dirac
deltas with equal weights.
\end{teo}

\proof
By Lemma~\ref{l_7.1+2}, the conical limit set
has full measure in $\partial \M H^k$ and $\Gamma$ acts ergodically on
$\partial \M H^k$.
Let $E$ be the set of points $z\in\partial \M H^k$ such that $\lambda_z$ is
the sum of two deltas with equal weights. We have to prove that $\mu_O(E)=0$.
Clearly, the set $E$ is $\Gamma$-invariant, and by
Lemma~\ref{l_5.0_7.3} it is measurable. Therefore it has either zero
or full $\mu_O$-measure.

Suppose that $E$ has full $\mu_O$-measure. Then, a $\mu_O$-measurable
map is well defined by 

$$\begin{array}{c}
f:\partial \M H^k\to\big(\partial\M H^n\times\partial\M
H^n\big)/\mfk S_2\\
z\mapsto \textrm{support of }\lambda_z.
\end{array}
$$
By Lusin theorem (see for example~\cite{afp:libro}) for all $\e>0$
there exists a compact set $A\subset\partial \M H^k$, with
$\mu_O(\partial \M H^k\setminus A)<\e$ and such that the restriction
of $f$ to $A$ is continuous. By Lemma~\ref{l_5.0_7.5}, for any
density-point $x$ of $A$, for $\mu_O\times\mu_O$-almost all
$(z_1,z_2)\in\partial\M H^k\times\partial\M H^k$ there exists a
sequence $\{\gamma_i\}\subset\Gamma$ such that for $j=1,2$
$$\gamma_i(z_j)\in A\qquad\textrm{ and }\qquad \gamma_i(z_j)\to x.$$
Therefore, for $j=1,2$
$$f(\gamma_i(z_j))\to f(x).$$
Up to pass to sub-sequences, $\rho(\gamma_i)$ converges either to an
isometry of $\M H^n$ or to a quasi-constant $C_a^b$
$\frac{~~~~}{~~~~}$
a quasi-constant is a map such that $C_a^b(p)=b$
for all points $p\neq a$ of $\wb{\M H}^n$
$\frac{~~~~}{~~~~}$
where the convergence is uniform
on compact sets not containing $a$ (see for example~\cite{Kap:libro}).

Let $\{\theta,\omega\}$ be the support of $\lambda_x$, that is
$\lambda_x=\frac{\delta_\theta+\delta_\omega}{2}$, and for $j=1,2$ let
$\{\xi_j,\zeta_j\}$ be the support of $\lambda_{z_j}$.

If $\rho(\gamma_i)$ converges to an isometry $\psi$, then for $j=1,2$
$$f(\gamma_i(z_j))=\rho(\gamma_i)\{\xi_j,\zeta_j\}\to\{\psi(\xi_j),\psi(\zeta_j)\}.$$
Since $f(\gamma_i(z_j))\to f(x)=\{\theta,\omega\}$, we get $f(z_1)=f(z_2)$.
On the other hand, if $\rho(\gamma_i)\to C_a^b$ and if $\xi_j\neq
a\neq\zeta_j$, we get
$$f(\gamma_i(z_j))\to\{b,b\}\neq f(x).$$

We have so proved that for $\mu_O\times\mu_O$-almost all
$(z_1,z_2)\in\partial\M H^k\times\partial\M H^k$ the supports of
$\lambda_{\z_1}$ and $\lambda_{z_2}$ share at least one point. Whence,
using Fubini's theorem and Lemma~\ref{l_5.0_7.5}, it
follows that there exists $\zeta\in\partial \M H^n$ such that for
$\mu_O$-almost all $z\in\partial\M H^k$, the support of $\lambda_z$
contains $\zeta$.
For $z\in\partial\M H^k$, let $\xi(z)$ denote the other point of the
support of $\lambda_z$, that is
$$f(z)=\{\zeta,\xi(z)\}.$$

For each $\gamma\in\Gamma$, for $\mu_O$-almost all $z$ we have
$f(z)=\{\zeta,\xi(z)\}$ and
$$\{\rho(\gamma)\zeta,\rho(\gamma)\xi(z)\}=\rho(\gamma)(f(z))=f(\gamma
z)=\{\zeta,\xi(\gamma z)\}.$$

Whence we have
\begin{itemize}
\item[Either:] the set $\{\zeta\}$ is $\rho(\Gamma)$-invariant.
\item[Or:] there exists $\gamma\in\Gamma$ such that
  $\rho(\gamma)\zeta\neq\zeta$, which implies that
  $\xi(z)=\rho(\gamma)^{-1}(\zeta)$ does not depend on $z$. Therefore
  the set $f(z)$ does not depend on $z$ and it is $\rho(\Gamma)$-invariant.
\end{itemize}
In both cases the image of $\rho$ is elementary, which contradicts the
hypotheses. It follows that the set $E$  cannot have full
$\mu_O$-measure.
\qed

\

Theorems~\ref{c_5.1_1.9} and~\ref{t_5.1_7.6} now easily follow.

\
\\
\noindent{\em Proof of Theorem~\ref{c_5.1_1.9}.}
It is well-known (see for example~\cite{Kap:libro}) that if $\M
H^k/\Gamma$ is a complete hyperbolic manifold of finite volume, then
$\delta(\Gamma)=k-1$ and $\Lambda_c$ has full-measure. Then, by
Lemma~\ref{l_7.1+2}, the hypotheses of
Theorem~\ref{t_lemmaF3} are satisfied, and the claim follows from Theorem~\ref{t_4}.\qed

\ \\
\noindent{\em Proof of Theorem~\ref{t_5.1_7.6}.}
The first claim is an immediate corollary of
Theorems~\ref{t_4} and~\ref{t_lemmaF3}.
Now, let $\f=\wb F\circ\wb G:\partial\M H^k\to\partial\M H^k$. The map
\f\ is clearly $\Gamma$-equivariant and $\mu_O$-measurable. Since the
identity is also $\Gamma$-equivariant and measurable, by
uniqueness it follows that $\wb F\circ\wb
G=\Id$ $\mu_0$-almost everywhere. The same holds replacing $\wb
F\circ\wb G$ with  $\wb G\circ\wb F$.
\qed

\

In particular this extends (and provides a new proof of) the following result.

\begin{teo}[\cite{cathu89,min94,som95}]\label{c_5.1_1.10}
 Let $M$ be a compact hyperbolic $3$-manifold fibred over the circle
  with fibre a surface $F$. Let $\pi_1(F)\simeq\Gamma<{\rm Isom}(\M
  H^2)$ and let $\rho:\Gamma\to{\rm Isom}(\M H^3)$ be the
  representation induced by the inclusion
  $\pi_1(F)\hookrightarrow\pi_1(M)$. Then, there exist measurable maps
  $\wb F:\partial \M H^2\to\partial\M H^3$ and  $\wb G:\partial \M
  H^3\to\partial\M H^2$ which are respectively $\rho$ and
  $\rho^{-1}$ equivariant. Moreover, almost everywhere
$$\wb F\circ\wb G={\rm Id}_{\M H^2}\qquad \wb G\circ\wb F={\rm Id}_{\M H^3}.$$
\end{teo}

\proof
Clearly, $\Gamma$ satisfies the hypotheses of Theorem~\ref{t_5.1_7.6}.
Moreover, by Lemma~\ref{l_7.1+2}
and~\cite[Corollary 2]{can93},
also $\pi_1(F)<\pi_1(M)<{\rm Isom}(\M H^3)$ satisfies the hypotheses of
Theorem~\ref{t_5.1_7.6}.\qed

\section{Convergence of Cannon-Thurston maps}

We begin this section by describing an example
of a  converging sequence of representations whose
Cannon-Thurston maps have no Cauchy sub-sequence with respect to the uniform
convergence. Then, we prove Theorem~\ref{t_7.0_1.8} (point-wise
 convergence almost everywhere.)

\begin{examp}[Souto]\label{ex_7.2_9.1}
\begin{rm}
Let $\Gamma$ be the fundamental group of a compact hyperbolic
surface. Then, there exists a sequence of discrete and faithful
representations $\rho_n:\Gamma\to\rm{Isom}(\M H^3)$ such that, if
$f_n$ and $f$ denote the corresponding Cannon-Thurston maps, then:
\begin{itemize}
\item $\rho_n(\Gamma)$ is quasi-Fuchsian
\item $\rho_n\to\rho$
\item $\rho(\Gamma)$ is geometrically finite
\item There is a converging  sequence of points $x_n\to x\in\partial\M
  H^2$ such that $f_n(x_n)\to y\neq f(x)$. In particular, no
  sub-sequence of $\{f_n\}$ converges uniformly to $f$.
\end{itemize}

To see that, let $AH$ denote set of discrete and faithful
representations $\rho:\Gamma\to\rm{Isom}(\M H^3)$ with the topology
of algebraic convergence. Let $\rho$ be a geometrically finite
representation in the closure of $AH$ (whence $\rho$ has accidental
parabolics.) The domain of discontinuity of $\rho(\Gamma)$ is
therefore non-empty, and we can thus pick a point $y$ in it. Since
the doubly degenerate representations are dense in the closure of
$AH$, there exists a sequence $\psi_n\to \rho$ with the property
that the limit set of $\psi_n(\Gamma)$ is the whole sphere $\partial
\M H^3$. Let $f_n$ and $f$ be the Cannon-Thurston maps corresponding
to $\psi_n$ and $\rho$. By equivariance, there exists a sequence of
points $x_n\in\Lambda(\Gamma)$ such that $x_n\to x\in
\Lambda(\Gamma)$ and $f_n(x_n)\to y$. Since $y\notin
\Lambda(\rho(\Gamma))$ we cannot have $f_n(x_n)\to f(x)$. Finally,
one can approximate the $\psi_n$'s with quasi-Fuchsian
representations $\rho_n$ getting the requested properties.
 \end{rm}
\end{examp}

\proof[Proof of Theorem~\ref{t_7.0_1.8}] For $z\in\Lambda(\Gamma)$
let $\delta_{f_i(z)}$ be the Dirac measure concentrated on $f_i(z)$.
For $x\in\M H^k$ consider the measures on $\partial \M
H^k\times\partial \M H^n$ defined by
$$\eta_{x,i}=\mu_x\times\{\delta_{f_i(z)}\}.$$
Up possibly to pass to sub-sequences, the measure $\eta_{O,i}$ have a
weak limit $\eta_O$, which we disintegrate as
$$\eta_O=\mu_O\times\{\beta_z\}$$
via some family of probability measures
$\{\beta_z\}_{z\in\Lambda(\Gamma)}$. Let now
$$\eta_x=\mu_x\times\{\beta_z\}.$$

It is readily checked that $\eta_{x,i}$ weakly converges to $\eta_x$
(because
$d\mu_x(\theta)=e^{-\delta(\Gamma)B_O(x,\theta)}d\mu_O(\theta)$.) It
follows that the measures $\eta_x$ are $\rho$-equivariant. Indeed,
for any fixed $\gamma\in\Gamma$ we have $\eta_{\gamma
x,i}=(\gamma,\rho_i(\gamma))_*\eta_{x,i}$, and $\rho_i(\gamma)$
converges to $\rho(\gamma)$ uniformly on $\partial \M H^n$.

As in Lemma~\ref{l4_4.3} this implies that the family $\beta_z$ is
$\rho$-equivariant, and uniqueness part of Theorem~\ref{t_4} implies
that in fact $\beta_z=\delta_{f(z)}$.

Now that we know that $\eta_x=\mu_x\times\{\delta_{f(z)}\}$, the
proof is completed by the fact that the weak convergence of
$\eta_{x,i}$ to $\eta_x$ is equivalent to the convergence of $f_i$
to $f$ almost everywhere (see for
example~\cite[Lemma~2.3]{alspre05}.) \qed

\end{document}